\documentclass{ifacconf}

\usepackage{graphicx}      
\usepackage{subcaption}
\usepackage{amsmath,amssymb,amsfonts}
\usepackage{natbib}        
\usepackage{bm, booktabs}
\newtheorem{remark}{Remark}
\newtheorem{assumption}{Assumption}
\usepackage{xcolor}
\begin{document}
\begin{frontmatter}

\title{Model-free source seeking of exponentially convergent unicycle: theoretical and robotic experimental results
  } 


\author[First]{Rohan Palanikumar} 
\author[First]{Ahmed A. Elgohary} 
\author[Second]{Victoria Grushkovskaya}
\author[First]{Sameh A. Eisa}

\address[First]{Department of Aerospace Engineering and Engineering Mechanics,
University of Cincinnati, Cincinnati, OH 45221 USA (palanire@mail.uc.edu, elgohaam@mail.uc.edu, eisash@ucmail.uc.edu).}
\address[Second]{Department of Mathematics,
        University of Klagenfurt, 9020 Klagenfurt am W\"orthersee, Austria (viktoriia.grushkovska@aau.at).}

\begin{abstract}
This paper introduces a novel model-free, real-time unicycle-based source seeking design. This design autonomously steers the unicycle dynamic system towards the extremum point of an objective function or physical/scalar signal that is unknown expression-wise, but accessible via measurements. A key contribution of this paper is that the introduced design converges exponentially to the extremum point of objective functions (or scalar signals) that behave locally like a higher-degree power function (e.g., fourth-degree polynomial function) as opposed to locally quadratic objective functions, the usual case in literature. We provide theoretical results and design characterization, supported by a variety of simulation results that demonstrate the robustness of the proposed design, including cases with different initial conditions and measurement delays/noise. Also, for the first time in the literature, we provide experimental robotic results that demonstrate the effectiveness of the proposed design and its exponential convergence ability. These experimental results confirm that the proposed exponentially convergent extremum seeking design can be practically realized on a physical robotic platform under real-world sensing and actuation constraints.
\end{abstract}

\begin{keyword}
Source Seeking, Extremum Seeking, Unicycle, Exponential Convergence, Robotic Experiment, Model-free Optimization, Light source seeking.
\end{keyword}

\end{frontmatter}

\section{Introduction}
Extremum seeking control (ESC) techniques are model-free, real-time control methods that drive a system to the optimum of a given objective function (\cite{ariyur2003real,review,guay2015time,scheinker2024100}), which may be unknown, provided it can be measured. ESC schemes are particularly attractive across various fields due to their minimal information requirements: they rely solely on applying perturbations to system parameters or inputs, and measuring the corresponding outputs of the objective function. Using these measurements in a feedback loop, ESC algorithms iteratively adjust the system input or parameter values to drive it towards the extremum, as can be seen in many applications, e.g., \cite{krstic2008extremum,dochain2011extremum,calli2012comparison,ECC_2024,elgohary2025model,elgohary2025extremum,elgohary2025hovering,pokhrel2022novel,eisa2023,grushkovskaya2018family,grushkovskaya2024design}.

One conventional problem that has increasingly involved the use of ESC methods is source-seeking (\cite{bajpai2024investigating,ghods2011extremum,zhu2013cooperative,li2014cooperative,bulgur2018light}). Source-seeking is the problem of steering a system (e.g., a mobile robot) autonomously to the maximum or minimum intensity of a physical/scalar signal present in a given domain using only on-board sensor measurements (e.g., heat, light, chemical concentration, among others). The advantage of ESC lies in the fact that it does not require knowledge of the signal's distribution or explicit modeling of the system to drive the system, autonomously, to the extremum point. Due to the model-free nature of ESC methods, source-seeking techniques based on ESC utilized a simple, basic unicycle model to represent the system within the control design (\cite{pokhrel2024extremum,matveev2011navigation,khong2014multi,cochran20093,grushkovskaya2018family,suttnerkrisitic2023unicycletorque,yilmaz2025unbiased,todorovski2024newton,elgohary2025model}). As a result, unicycle-based methods for source-seeking often resort to Lie bracket averaging (\cite{durr2013lie,scheinker2017bounded,grushkovskaya2018class,ghadiri2020normalized,pokhrel2023higher}) in theoretical analysis and design. It is worth noting that the aforementioned unicycle-based designs for source-seeking are based on ESC laws that are gradient-based, which is more effective (e.g., guarantee exponential convergence) if the objective function or the scalar signal behaves locally like a quadratic function.    

On the other hand, optimization problems involving fourth-order polynomial functions have become an increasingly prevalent topic in the literature (\cite{qi2004global, zhang2015optimality}). In fact, some emerging applications in science and engineering involve the need to optimize and control a fourth-order polynomial function, such as but not limited to the distribution of the coefficient of performance for power optimization ($C_p$) in wind turbines (\cite{eisa2019modeling}) and pure-quartic solitons in optics, which exhibit fourth-order dispersion and have possible applications with lasers (\cite{soltani2023pure}).



\textbf{Motivation and Contribution.} 
In our recent work \cite{grushkovskaya2025extremum}, new ESC laws were derived based on higher-order Lie bracket averaging. This paved the way for the prospect of having exponentially convergent ESC laws for objective functions that behave locally not as quadratic functions, but as power functions with higher-order degrees. Inspired by the prospect of realizing exponential convergence for high-order objective functions, we seek to utilize the ESC method in \cite{grushkovskaya2025extremum} and propose a novel, first-of-its-kind, exponentially convergent unicycle-based design for source-seeking. The main contribution of this work is the development and experimental validation of an exponentially convergent unicycle-based ESC law that realizes higher-order Lie bracket motion in real-time. In this preliminary study, we only focus on source-seeking of objective functions that behave locally as a fourth-degree power function. We prove stability of the proposed design and provide a variety of simulation results that demonstrate the robustness of the proposed design, including cases with different initial conditions and measurement delays/noise. Additionally, we provide experimental results using a TurtleBot3 robot (similar to \cite{ECC_2024,elgohary2025model}) and compare the proposed third-order Lie bracket-based design with traditional first-order Lie bracket-based designs. We also perform a light source-seeking experiment in a complete model-free fashion. Our results illustrate the potential and effectiveness of the proposed exponentially convergent unicycle design for source-seeking.

\section{Brief Background}
In this section, we provide brief context regarding the exponentially convergent ESC law that we use for higher-order objective functions. We consider the two-input control-affine ESC of the following form, as introduced in \cite{grushkovskaya2025extremum}:
\begin{equation}
   \dot{x} = g_1(J(x))u_1^\varepsilon(t) + g_2(J(x))u_2^\varepsilon(t),
   \label{eq:two_inp_esc_system}
\end{equation}
where $x \in D \subset \mathbb{R}^n$ denotes the system's state, $J : D \to \mathbb{R}$ is the objective function, 
$u_1^\varepsilon(t) = \varepsilon^{1/N -1}v_1(t/\varepsilon)$, 
$u_2^\varepsilon(t) = \varepsilon^{1/N -1}v_2(t/\varepsilon)$, 
$N-1 \in \mathbb{N}$ denotes the order of the Lie bracket to be excited, and $\varepsilon>0$ is a small positive parameter, such that the dithers $v_1(t/\varepsilon)$ and $v_2(t/\varepsilon)$ excite the Lie bracket
\[
g_{I_N}(z)=\big[\big[\dots[g_1,\underbrace{g_2],g_2\big],\dots,g_2\big]}_{N-1\ \text{times}}\dots]\big](z)
\]
at time $t = \varepsilon$ with $I_N = (1,\underbrace{2,\dots,2}_{N-1\ \text{times}})$, where
\begin{equation}
[g_i,g_j]:=\frac{\partial g_j}{\partial x}g_i-\frac{\partial g_i}{\partial x}g_j.\nonumber
\end{equation}
Moreover, we adopt all mathematical assumptions on the vector fields $g_1$, $g_2$ provided in \cite[Section II]{grushkovskaya2025extremum} and recall that $g_1,g_2:\mathbb{R}\to\mathbb{R}^n$ are smooth vector fields evaluated at the scalar argument $J(x)$.
However, we state below the assumption clarifying the structure and class of objective functions admitted in \cite{grushkovskaya2025extremum} as well as in this paper:

\begin{assumption}
The function $J\in C^m(D,\mathbb{R})$, where $D\subset\mathbb{R}^n$ and $m\geq 2$, then:
\begin{itemize}
    \item[A1.1)] there exists an $x^*{\in} D$  such that 
    $$\nabla J(x)=0 \text{ if and only if }x=x^*,$$
    and $J(x^*)=J^*{\in}\mathbb R$, $J(x)>J(x^*)$ for all $x\in D{\setminus}\{x^*\}$.
    \item[A1.2)] There exist constants $\alpha_1,\alpha_2,\beta_1,\beta_2,\mu$ such that, for all {$x\in D$},
     $$ \begin{aligned}
     \alpha_1\|x{-}x^*\|^{m} \le &J(x)-J^* \le \alpha_2\|x{-}x^* \|^{m},\\
      \beta_1 (J(x)-J^*)^{1{-}\frac{1}{m}}\le  & \|\nabla J(x)\|\le\beta_2 (J(x)-J^*)^{1{-}\frac{1}{m}},\\
    \left\|\frac{\partial^2 J(x)}{\partial x^2}\right\|{\le}&\mu (J(x)-J^*)^{1{-}\frac{2}{m}}.
    \end{aligned}
    $$
    Here, $\|\cdot\|$ denotes the Euclidean norm.
\end{itemize}
\end{assumption}

\begin{remark}
    Assumption A1.1 states that $J$ has an isolated minimum, $x^*$, at which it has a value $J^*$. Assumption A1.2 provides the requirement that the admitted objective function, $J$, behaves locally as a power function of degree $m$; this is consistent with ESC literature (e.g., \cite{grushkovskaya2018class}).
\end{remark}

The selection of the dither signals is done as in \cite{grushkovskaya2024design,grushkovskaya2025extremum}, conditions C1--C3:
\begin{itemize}
    \item \textbf{C1.} \textcolor{black}{With} $v_1(t/\varepsilon)=2\sqrt{\kappa_{12}\pi}\cos{(2\kappa_{12}\pi t/\varepsilon)}$ and $v_2(t/\varepsilon)=2\sqrt{\kappa_{12}\pi}\sin{(2\kappa_{12}\pi t/\varepsilon)}$, $\kappa_{12} \in \mathbb{Z}$, the first-order Lie bracket is excited.
    \item \textbf{C2.} \textcolor{black}{With} $v_1(t/\varepsilon)=-2(4\kappa_{122}\pi)^{2/3}\cos{(4\kappa_{122}\pi t/\varepsilon)}$ and $v_2(t/\varepsilon)=(4\kappa_{122}\pi)^{2/3}\cos{(2\kappa_{122}\pi t/\varepsilon)}$, $\kappa_{122} \in \mathbb{Z}$, the second-order Lie bracket is excited.
    \item \textbf{C3.} \textcolor{black}{With} $v_1(t/\varepsilon)=6(2\kappa_{1222}\pi)^{3/4}\sin{(6\kappa_{1222}\pi t/\varepsilon)}$ and $v_2(t/\varepsilon)=2(2\kappa_{1222}\pi)^{3/4}\cos{(2\kappa_{1222}\pi t/\varepsilon)}$, $\kappa_{1222} \in \mathbb{Z}$, the third-order Lie bracket is excited.
\end{itemize}
From \cite{grushkovskaya2025extremum}, $g_1,\ g_2$ are chosen according to:
\begin{equation}
    g_{I_N}(J(x))=-c_N J^{(N-1)}(x),
    \label{eq:req_g_IN}
\end{equation}
with the simple selection to satisfy \eqref{eq:req_g_IN} being $g_1(z) = -1^{(N+1)}z$ and $g_2(z) = 1$. It is also to be noted that the conditions C1--C3 above guarantee that only the desired Lie bracket is excited while all others vanish per \cite{Gau14,pokhrel2023higher,grushkovskaya2024design}. In this paper, we are interested in third-order Lie bracket excitation (i.e., condition C3).

\section{Main Results: Proposed Design and Stability Analysis}
We now propose our exponentially convergent unicycle design and prove its stability.
\subsection{Proposed Design}
\begin{figure}[t]
    \centering
    \includegraphics[width=1\linewidth]{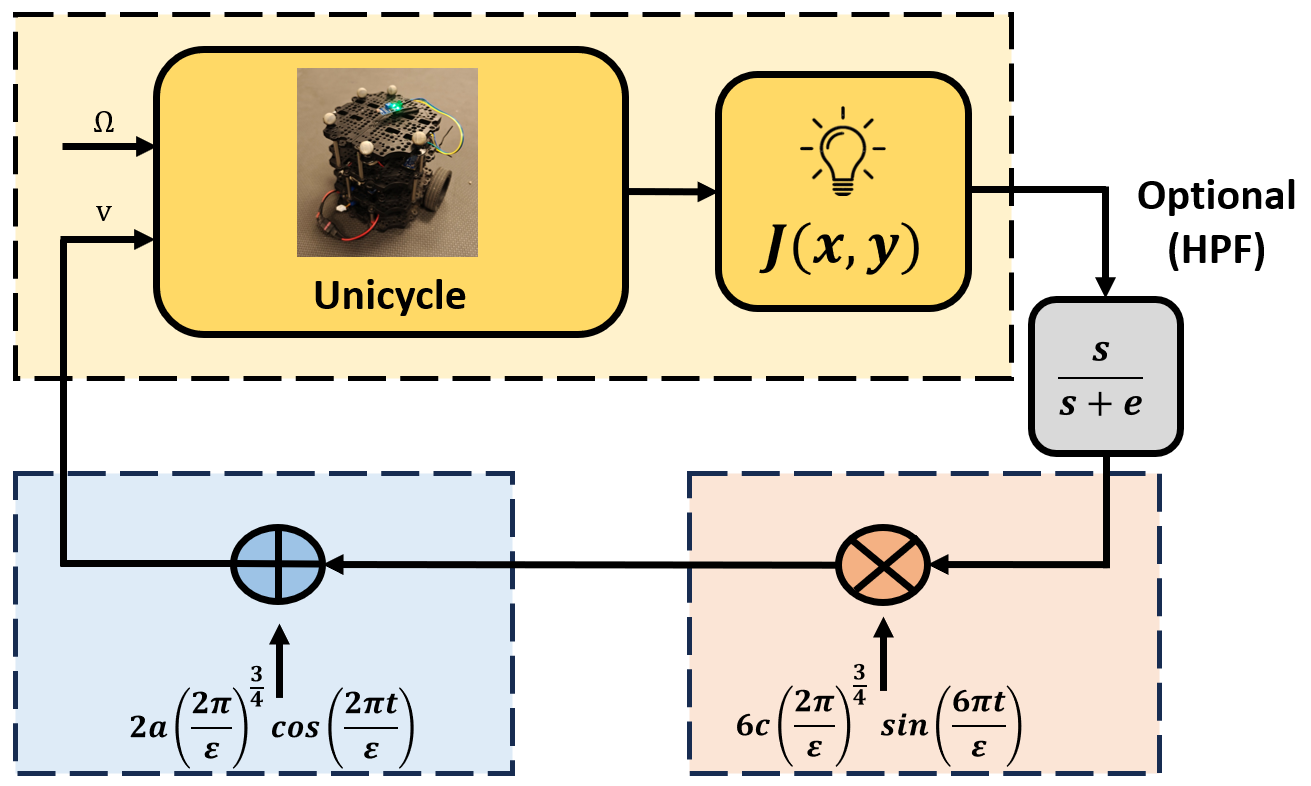}
    \caption{The proposed ESC design for exponentially convergent unicycle. We use the proposed unicycle design for differential drive robotic experiments. In experiments and simulations we use mathematically known objective function and unknown light source.}
    \label{fig:ESC_diagram}
\end{figure}

We \textcolor{black}{aim to} use an ESC law based on condition C3, provided in the previous section, to achieve an exponentially convergent unicycle design. Let us consider the kinematic differential equations for unicycle dynamics with constant angular velocity as follows (\cite{elgohary2025model,grushkovskaya2018family}):
\begin{equation}\label{eq:main_unicycle}
\begin{aligned}
    \dot{x} &= \mathrm{v}\cos(\Omega t), \\
    \dot{y} &= \mathrm{v}\sin(\Omega t),
\end{aligned}
\end{equation}
where $x$ and $y$ describe the current position/coordinates, $\mathrm{v}$ is the linear/\textcolor{black}{translational} velocity, and $\Omega$ is the angular velocity.  
While the results can be extended to any order, in this study, we only focus \textcolor{black}{on} the unicycle design to address objective functions that behave locally as \textcolor{black}{a} fourth\textcolor{red}{-}order degree polynomial in consistency with \cite[Assumption 1]{grushkovskaya2025extremum}. Hence, we suppose:
\begin{equation}
    J(x,y) = C_1(x-x_d)^4 + C_2(y-y_d)^4,
    \label{eq:objective_function}
\end{equation}
where $C_1,C_2>0$. As analyzed, shown and simulated in \cite{grushkovskaya2025extremum}, we need to excite \textcolor{black}{the} third-order Lie bracket for exponential convergence since the objective function order is \textcolor{black}{of} a fourth-degree. Hence, we choose our control law for $\mathrm{v}$ of the unicycle based on condition C3. That is,
\begin{equation}
    \mathrm{v} = 2(2\pi / \varepsilon)^{(3/4)} (3cJ(x,y) \sin{(6\pi t/\varepsilon)} + a\cos{(2\pi t/\varepsilon)}),
    \label{eq:linear_vel_input}
\end{equation}
where $\varepsilon>0$ is a small parameter and $a,c>0$ are design constants.
Now we are in a position to put all elements of the proposed design together, which is depicted in Figure \ref{fig:ESC_diagram}. In state space representation, including an optional high-pass filter (HPF) into the system (\cite{ECC_2024,elgohary2025model}), the proposed design becomes:
 \small
 \begin{equation}
 \begin{aligned}
     \dot{x} &= \big(2(2\pi / \varepsilon)^{(3/4)} (3c(J-eh) \sin{(6\pi t/\varepsilon)} + a\cos{(2\pi t/\varepsilon)})\big)\cos{(\Omega t)}, \\
     \dot{y} &= \big(2(2\pi / \varepsilon)^{(3/4)} (3c(J-eh) \sin{(6\pi t/\varepsilon)} + a\cos{(2\pi t/\varepsilon)})\big)\sin{(\Omega t)}, \\
     \dot{h} &= J(x,y) - eh,
 \end{aligned}
 \label{eq:full_unicycle_system_with_hpf}
 \end{equation}
 \normalsize
 where $h$ is the optional filter state and $e$ is the filter constant. If the design is to be considered without the HPF, we set $h=0$ and omit its dynamical equation. 

\subsection{Stability Analysis}
In this subsection, we provide the stability analysis for the proposed unicycle design \eqref{eq:full_unicycle_system_with_hpf}. We will follow the traditional methodology in Lie bracket-based ESC literature (e.g., \cite{durr2013lie,grushkovskaya2018class,pokhrel2023higher}) where the stability property of the ESC system is characterized by the corresponding Lie bracket system (LBS). The corresponding LBS to \eqref{eq:full_unicycle_system_with_hpf} without the optional HPF (i.e., $h=0$) is a third-order LBS based on condition C3 of the form (\cite{grushkovskaya2025extremum}):
$$
\dot{\bar{x}}=\big[[[g_1,g_2],g_2],g_2\big]
$$
with
$$g_1=\left(
\begin{array}{c}
  cJ(\bar{x},\bar{y})\cos(\Omega t) \\
 cJ(\bar{x},\bar{y})\sin(\Omega t) 
\end{array}
\right), \,
g_2=\left(
\begin{array}{c}
  a\cos(\Omega t) \\
  a\sin(\Omega t)
\end{array}
\right).$$
Let $J_{\bar{x}\bar{x}\bar{x}}$ denotes the third order partial derivative of $J(\bar{x},\bar{y})$ with respect to $\bar{x}$ and $J_{\bar{y}\bar{y}\bar{y}}$ denotes the third order partial derivative of $J(\bar{x},\bar{y})$ with respect to $\bar{y}$. Then,
\begin{equation}\label{LBS}
\begin{aligned}
\dot{\bar{x}}&=-ca^3\cos(\Omega t)(J_{\bar{x}\bar{x}\bar{x}}\cos^3(\Omega t)+J_{\bar{y}\bar{y}\bar{y}}\sin^3(\Omega t))\\
&=-\cos(\Omega t)(c_1(\bar{x}-\bar{x}_d)\cos^3(\Omega t)+c_2(\bar{y}-\bar{y}_d)\sin^3(\Omega t)),\\
\dot{\bar{y}}&=-ca^3\sin(\Omega t)(J_{\bar{x}\bar{x}\bar{x}}\cos^3(\Omega t)+J_{\bar{y}\bar{y}\bar{y}}\sin^3(\Omega t))\\
&=-\sin(\Omega t)(c_1(\bar{x}-\bar{x}_d)\cos^3(\Omega t)+c_2(\bar{y}-\bar{y}_d)\sin^3(\Omega t)),
\end{aligned}
\end{equation}
with $c_1=4!ca^3C_1$, $c_2=4!ca^3C_2$. 
The following result establishes the stability properties of system~\eqref{LBS}.
\begin{thm}\label{thm1}
\textit{Let   one of the following conditions be satisfied: 
\begin{itemize}
    \item [i)] $c_1\in(3c_2/5,c_2]$, $\Omega>\frac{2c_1(3c_2-c_1)(2c_2-c_1)}{16c_1^2-(3c_2-c_1)^2}$;
    \item [ii)] $c_2\in(3c_1/5,c_1]$, $\Omega>\frac{2c_2(3c_1-c_2)(2c_1-c_2)}{16c_2^2-(3c_1-c_2)^2}$.
\end{itemize}
Then the equilibrium $x^*=x_d,y^*=y_d$ of system~\eqref{LBS} is exponentially stable. }
\end{thm}

\textbf{Proof.} We begin with the following change of variables:
$$
\begin{aligned}
 & \xi=(\bar{x}-\bar{x}_d)\cos(\Omega t)+(\bar{y}-\bar{y}_d)\sin(\Omega t),\\
 & \eta=(\bar{x}-\bar{x}_d)\sin(\Omega t)-(\bar{y}-\bar{y}_d)\cos(\Omega t),
\end{aligned}
$$
with the inverse transformation
$$
\begin{aligned}
 & \bar{x}-\bar{x}_d=\xi\cos(\Omega t)+\eta\sin(\Omega t),\\
 & \bar{y}-\bar{y}_d=\xi\sin(\Omega t)-\eta\cos(\Omega t).
\end{aligned}
$$
Observe that 
$$
\begin{aligned}
 \dot\xi=&\dot{\bar{x}}\cos(\Omega t)+\dot{\bar{y}}\sin(\Omega t)\\
 &-\Omega((\bar{x}-\bar{x}_d)\sin(\Omega t)-(\bar{y}-\bar{y}_d)\cos(\Omega t)) \\
 =&-(c_1(x-x_d)\cos^3(\Omega t)+c_2(y-y_d)\sin^3(\Omega t))-\Omega \eta\\
 =&-\xi(c_1\cos^4(\Omega t)+c_2\sin^4(\Omega t))\\
 &-\frac{1}{2}\eta\sin(2\Omega t)(c_1\cos^2(\Omega t)-c_2\sin^2(\Omega t)),
\end{aligned}
$$
and
$$
\begin{aligned}
 \dot\eta=&\dot x\sin(\Omega t)-\dot y\cos(\Omega t)\\
 &+\Omega((x-x_d)\cos(\Omega t)+(y-y_d)\sin(\Omega t)) \\
 =&\Omega \xi.\\
\end{aligned}
$$
Thus, in the new variables system~\eqref{LBS} takes the form
\begin{equation}\label{LBS_xi}
  \begin{aligned}
    &\dot\xi=-\kappa_1(t)\xi-\eta(\Omega+\kappa_2(t)),\\
    &\dot \eta=\Omega \xi.
  \end{aligned}
\end{equation}
where 
$$
\kappa_1(t)=c_1\cos^4(\Omega t)+c_2\sin^4(\Omega t),
$$
$$
\kappa_2(t)=\frac{1}{2}\sin(2\Omega t)(c_1\cos^2(\Omega t)-c_2\sin^2(\Omega t)).
$$
Using the identities $$\cos^4(\Omega t)+\sin^4(\Omega t)=1-\frac{1}{2}\sin^2(2\Omega t)$$ and 
$$
c_1\cos^2(\Omega t)-c_2\sin^2(\Omega t)=\frac{1}{2}(c_1-c_2)+(c_1+c_2)\cos(2\Omega t),
$$
we obtain the following estimates for the coefficients of system~\eqref{LBS_xi}: for all $t\ge0$,
\begin{equation}
    \label{kappa_est}
\begin{aligned}
 k_{11}\le&\kappa_1(t)\le k_{12} ,\\
 |&\kappa_2(t)|\le k_2, 
\end{aligned}   
\end{equation}
where $k_{11}=\frac{1}{2}\min\{c_1,c_2\}$, $k_{12}=\max\{c_1,c_2\}$, $k_2=\frac{1}{4}|c_1-c_2|+\frac{1}{8}(c_1+c_2)$.
To prove the exponential stability of the trivial solution of system~\eqref{LBS_xi}, consider the function 
\begin{equation}\label{Lyap}
  V(\xi,\eta)=\frac{1}{2}\xi^2+\frac{1}{2}\eta^2+\gamma\xi\eta
\end{equation}
with $\gamma\in(0,1)$ to be defined. 
Then
$$
\begin{aligned}
  \dot V=&-(\kappa_1(t)-\gamma\Omega)\xi^2-\gamma(\Omega+\kappa_2(t))\eta^2\\
  &-(\gamma\kappa_1(t)+\kappa_2(t))\xi\eta\\
  \le&-\alpha_1\xi^2-\alpha_2\eta^2+\alpha_{12}|\xi\eta|\\
\end{aligned}
$$
with
$$
\alpha_1=(k_{11}-\gamma\Omega),\,\alpha_2=\gamma(\Omega+k_2),\,\alpha_{12}=(\gamma k_{12}+k_2).
$$
Requiring $\gamma<\frac{k_{11}}{\Omega}$, we ensure $\alpha_1>0$. Thus, to have negative definiteness of the function $\dot V$ it is enough to ensure 
$$
\alpha_{12}^2-4\alpha_1\alpha_2<0,
$$
which is equivalent to the requirement
\begin{equation}
    \label{ineqD1}
    (4\Omega^2+4k_2\Omega+k_{12}^2)\gamma^2-2(2\Omega k_{11}+2k_{11}k_2-k_2k_{12})\gamma+k_2^2<0.
\end{equation}
To ensure that the latter inequality is solvable, it is enough to have
$$
2(2\Omega k_{11}+2k_{11}k_2-k_2k_{12})^2-4k_2^2(4\Omega^2+4k_2\Omega+k_{12}^2)\ge 0.
$$
Factorizing the left hand side, we obtain
$$
16(\Omega+k_2)(\Omega(k_{11}^2-k_2^2)+k_{11}k_2(k_{11}-k_{12}))\ge0,
$$
which, in turn, leads to the requirement
\begin{equation}
    \label{ineqD2}
    \Omega(k_{11}^2-k_2^2)+k_{11}k_2(k_{11}-k_{12})\ge0.
\end{equation}
Because of the definition of $k_{11},k_{12}$, the difference $k_{11}-k_{12}$ is always negative, while by the conditions of the Theorem, 
$k_{11}^2-k_2^2>0$ and 
$$
\Omega\ge \frac{k_{11}k_2(k_{12}-k_{11})}{k_{11}^2-k_2^2}.
$$
Thus, inequality~\eqref{ineqD2} is satisfied, which means that there exists a $\hat \gamma>0$ such that, for all $\gamma\in (0,\hat\gamma]$, requirement~\eqref{ineqD1} is satisfied. Thus, if $\gamma\in(0,\min\{1,\hat\gamma,\frac{k_{11}}{\Omega}\}]$, then 
$-\alpha_1\xi^2-\alpha_2\eta^2+\alpha_{12}|\xi\eta|$ is negative definite, therefore,  such that 
$$\dot V\le -\mu (\xi^2+\eta^2),$$
where
$\mu >0$ is the greatest eigenvalue of the matrix $\left(\begin{matrix}
    \alpha_1 & \alpha_{12}/2\\
    \alpha_{12}/2 &\alpha_2
\end{matrix}\right)$. Similarly,
$$
\frac{1-\gamma}{2}(\xi^2+\eta^2)\le V\le \frac{1+\gamma}{2} (\xi^2+\eta^2).
$$
Thus,
$$
\dot V\le -\frac{2\mu}{1+\gamma}V,
$$
which yields the exponential decay
$$
V(t)\le V(0)e^{-\frac{2\mu t}{1+\gamma}},
$$
and
$$
\xi^2+\eta^2\le \frac{1+\gamma}{1-\gamma} (\xi(0)^2+\eta(0)^2)e^{-\frac{2\mu t}{1+\gamma}}.
$$
Coming back to the $(\bar{x},\bar{y})$-variables, we conclude
$$
(\bar{x}-\bar{x}_d)^2+(\bar{y}-\bar{y}_d)^2\le \frac{1+\gamma}{1-\gamma} (\xi(0)^2+\eta(0)^2)e^{-\frac{2\mu t}{1+\gamma}}.
$$
\begin{remark}
\textit{The simplest case in which the conditions of Theorem~1 are satisfied is when $c_1=c_2=c>0$, $\Omega>\frac{c}{3}$.
  We emphasize that the  conditions on $C_1,C_2,\Omega$ are sufficient but not necessary, as they result from the particular choice of the Lyapunov function used in the proof. We expect that more general conditions could be derived by using, for example, a Lyapunov function with time-periodic coefficients or by applying Barbalat’s lemma. We leave this question for future work. }
\end{remark}
\begin{thm}
\textit{The unicycle system \eqref{eq:full_unicycle_system_with_hpf} with $h=0$ is practically exponentially stable for any compact set $D\subset \mathbb{R}^2$ such that $(x_d,y_d) \in D$}.  
\end{thm}
\textbf{Proof.} Let us define the vectors $\bm{X}=(x,y)\in \mathbb{R}^2$, $\bm{X}_d=(x_d,y_d)\in \mathbb{R}^2$ and $\bm{\bar{X}}=(\bar{x},\bar{y})\in \mathbb{R}^2$. We have:
\begin{align}\label{eq:exp_pract_ineq}
    ||\bm{X}-\bm{X}_d||&=||\bm{X}-\bm{\bar{X}}+\bm{\bar{X}}-\bm{X}_d||\\
    &\leq ||\bm{\bar{X}}-\bm{X}_d||+||\bm{X}-\bm{\bar{X}}||. \nonumber
\end{align}
For any initial condition $\bm{X}_0=\bm{X}(0)=\bm{\bar{X}}(0)$, from Theorem 1, there exists $q_1>0$ and $q_2>0$ such that $||\bm{\bar{X}}-\bm{X}_d||\leq q_1e^{-q_2t}$. Moreover, per \cite[Theorem 4]{pokhrel2023higher}, there exists $d(\varepsilon)\to 0$ as $\varepsilon \to 0$ such that $||\bm{X}-\bm{\bar{X}}||\leq d(\varepsilon)$ for all $T>0$ and $t\in [0,T]$. Hence, the inequality \eqref{eq:exp_pract_ineq} becomes:
\begin{equation}\label{eq:practical_stability}
    ||\bm{X}-\bm{X}_d||\leq q_1e^{-q_2t} + d(\varepsilon), \: \: \forall T>0 \: \text{and} \: t\in [0,T]. 
\end{equation}
\begin{remark}
   \textit{ The inequality \eqref{eq:practical_stability} guarantees that for any positive time $T$, which can be made large as needed, the trajectories of the unicycle system \eqref{eq:full_unicycle_system_with_hpf} will decay exponentially to a neighborhood about the extremum $\bm{X}_d$. Said neighborhood can be made arbitrarily small via the parameter $\varepsilon$. The reader can refer to \cite{khalil2002,Maggia2020higherOrderAvg,pokhrel2023higher} for more details on the concept of practical exponential stability.}
\end{remark}

\subsection{Tuning Control Parameters}
\label{subsec: tuning_control_parms}
For tuning the exponentially convergent unicycle system, we provide some basic guidelines for the three control parameters that need to be tuned: $c$, $a$, and $\varepsilon$. The parameter $c$ acts as a learning rate parameter that is scaling the feedback of the objective function measurement and is affecting the rate of the velocity $v$. In simulation, better flexibility can be expected for $c$ as one can theoretically increase, as needed, the rate by which the feedback measurement of $J$ affects the rate of $v$. However, in real-world applications, there is a hardware limitation for how $v$ can change (i.e., limitation for the rate of $v$), hence $c$ can be increased carefully. The parameter $a$ is a scale of the added perturbation signal that is needed for Lie bracket realization. In simulations, balance between $a$ and $c$ can be considered. However, in real-world experiments, $a$ can be chosen such that it does not cause hardware violation for the rate of $v$, so $a$ can always be chosen after figuring out the viable value(s) of $c$. 

As far as $\varepsilon$ is concerned, in theory/simulations, the smaller the $\varepsilon$, the better; however, this is only possible in theory/simulations and is limited by the practical perturbation frequency the hardware implementation can admit. In particular, $\varepsilon$ determines the perturbation frequency $\omega=\frac{2\pi}{\varepsilon}$, and excessively small values of $\varepsilon$ correspond to high-frequency dither signals that may exceed actuator bandwidth, violate the Nyquist limit associated with the fixed control-loop sampling rate, or induce additional computational delays. Consequently, $\varepsilon$ is selected as the smallest value that can be reliably implemented while maintaining stable command tracking.

Furthermore, we note that the selection of $c$ and $a$ must be in accordance with Theorem~\ref{thm1}, such that either condition $\bm{i}$ or $\bm{ii}$ is satisfied. Lastly, when selecting $\varepsilon$, it must be sufficiently small and less than the threshold, $\varepsilon^{*}$, where $\varepsilon^{*}$ is not explicitly known in averaging analysis context (\cite{Maggia2020higherOrderAvg}). Additionally, in practice $\varepsilon$ may be selected by first identifying an admissible perturbation frequency $\omega$ based on hardware and sampling constraints, and then computing $\varepsilon=\frac{2\pi}{\omega}$.

\section{Simulation and Experimental Results}
\begin{figure}[t]
    \centering
    \includegraphics[width=\linewidth]{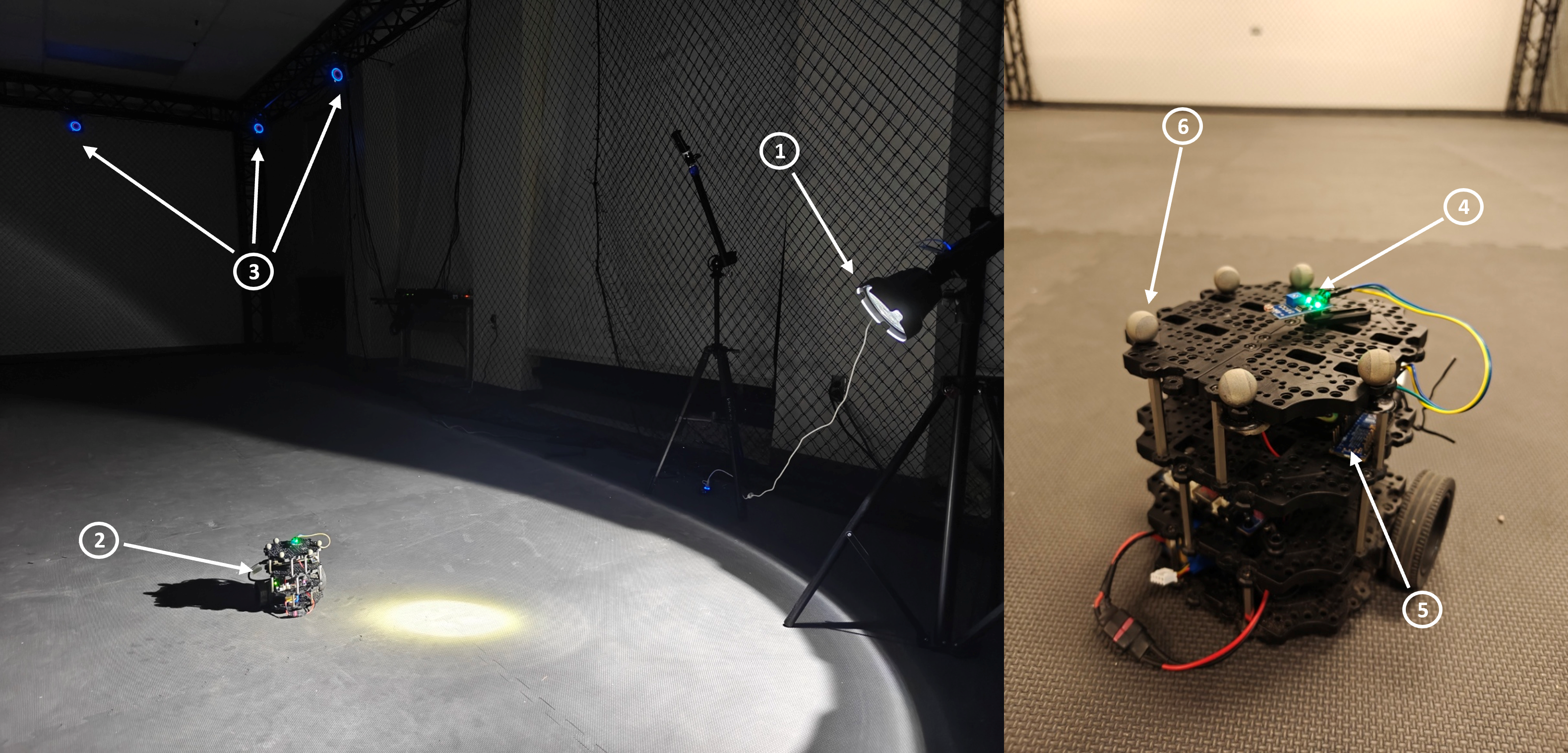}
    \caption{Modeling, Dynamics, and Control Lab (\cite{MDCL_Lab}): (1) light source; (2) TurtleBot3 robot; (3) motion capture system (MCS); (4) analog light sensor; (5) Arduino Nano ESP32; (6) MCS markers.}
    \label{fig:experimental_setup}
\end{figure}

In this section, we present the simulation results, the experimental setup, and the experimental validation for the proposed exponentially convergent unicycle design. To start, we briefly review the experimental setup, as shown in Figure~\ref{fig:experimental_setup}.

\subsection{Experimental Setup}
For experiments where the objective function is known, a motion capture system (MCS) is used to track the robot’s position in the $(x, y)$ coordinates. For the light source-seeking experiment, an analog light sensor connected to an Arduino Nano ESP32 board is mounted on the robot to measure light intensity. The robot used in our experiment is TurtleBot3 (see detailed information about this robot in \cite{bajpai2024investigating,ECC_2024}). The sensor readings are transmitted to the computer and used as feedback for the ESC system. The MCS is also used in this case to record and observe the robot’s position for performance evaluation. The light sensor consists of a small photoresistor that measures the intensity of light incident on its surface. As the light intensity increases, the sensor output decreases; hence, minimizing the sensor’s measured value corresponds to approaching the point of maximum light intensity.

\subsection{Simulation Results}
\begin{figure*}[t]
    \centering
    \begin{subfigure}[b]{0.33\linewidth}
        \centering
        \includegraphics[width=\linewidth]{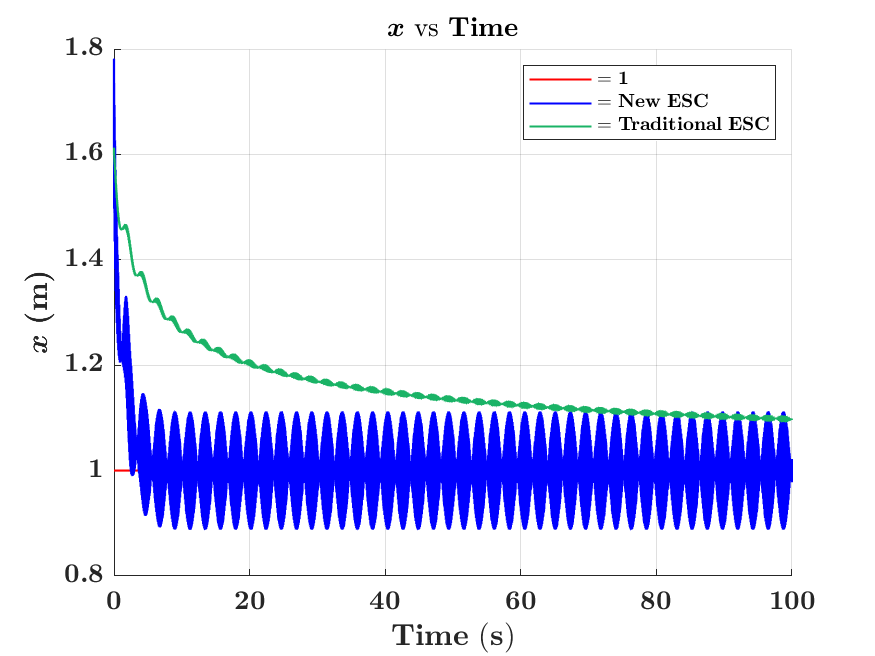}
        \caption{}
    \end{subfigure}\hfill
    \begin{subfigure}[b]{0.33\linewidth}
        \centering
        \includegraphics[width=\linewidth]{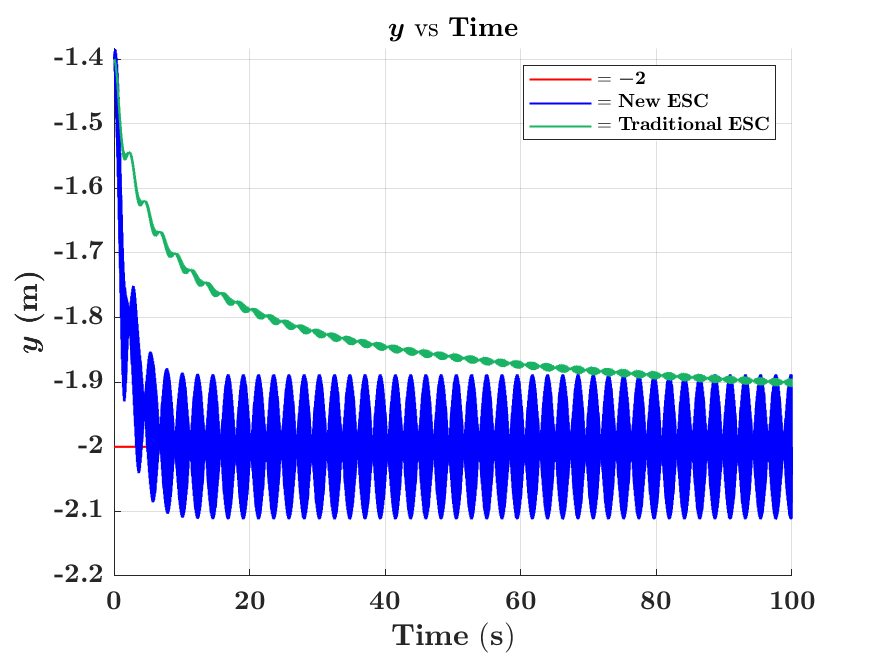}
        \caption{}
    \end{subfigure}\hfill
    \begin{subfigure}[b]{0.33\linewidth}
        \centering
        \includegraphics[width=\linewidth]{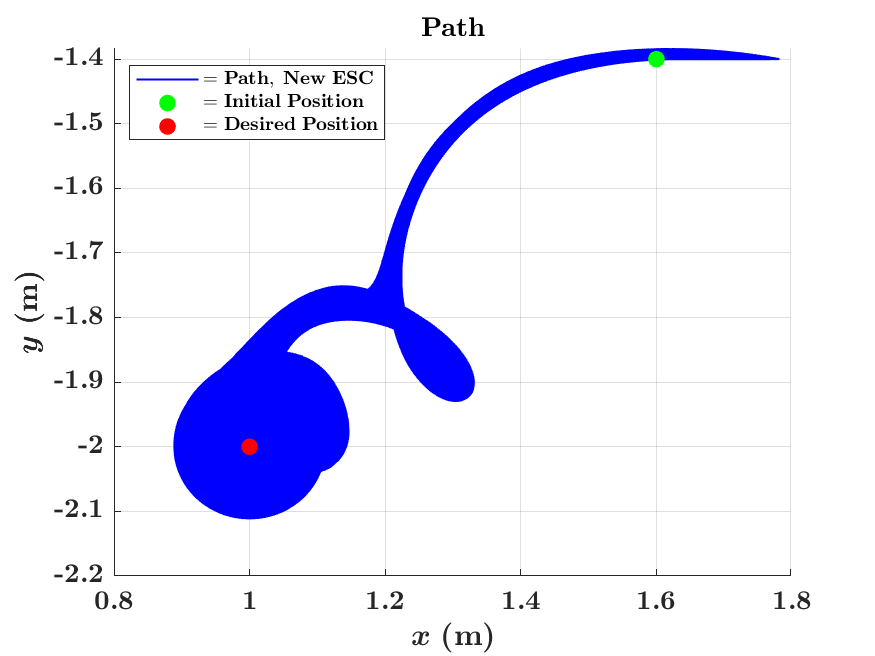}
        \caption{}
    \end{subfigure}
    \caption{Simulation results comparing the proposed exponentially convergent unicycle ESC with the traditional ESC based on first-order Lie bracket design from \cite{durr2013lie} for a fourth-order objective function: (a) $x$ position, (b) $y$ position, and (c) planar trajectory of the proposed unicycle design.}
    \label{fig:simulation_results}
\end{figure*}

In this subsection, we compare the traditional ESC design commonly found in the literature that is based on first-order Lie bracket design (\cite{durr2013lie}) with the proposed exponentially convergent unicycle ESC design developed for a fourth-order objective function based on third-order Lie bracket. The objective function used in this study is defined as
\[
J(x,y) = (x-1)^4 + (y+2)^4,
\]
which attains its minimum at $(x_d, y_d) = (1, -2)$. We remark here that this objective function satisfies the condition required for Theorem 1 (see Remark 1). The optional high-pass filter (HPF) is disabled in these simulations to isolate the effect of the proposed control law (i.e., $h=0$). The complete set of simulation parameters is provided in Table~\ref{tab:tabsimul}, and identical values are used for both ESC designs to ensure a fair comparison.
\begin{table}[t]
\centering
\caption{Simulation Parameters}\label{tab:tabsimul}
\begin{tabular}{@{}ll@{}}
\toprule
\textbf{Parameter} & \textbf{Value} \\ 
\midrule
$C_1, \ C_2$ & $1,\ 1$ \\
$a$ & $0.5$ \\
$c$ & $0.5$ \\
$\varepsilon$ & $0.001$ \\
$\Omega$ & $1.4\ \text{rad/s}$ \\
$x(0), \ y(0)$ & $1.6,\ -1.4\ \text{m}$ \\
$x_d, \ y_d$ & $1,\ -2\ \text{m}$ \\
\bottomrule
\end{tabular}
\end{table}
Figure~\ref{fig:simulation_results} illustrates the time histories of $x$ and $y$ as well as the planar trajectory of the unicycle. The results clearly demonstrate that the proposed exponentially convergent ESC achieves significantly faster convergence toward the desired equilibrium compared with the traditional ESC method from literature based on first-order Lie bracket design (\cite{durr2013lie}). Specifically, for the fourth-order objective function considered, the traditional ESC fails to converge within the 100~s simulation window, whereas the proposed design successfully converges within approximately 20~s. Lastly, the proposed exponentially convergent unicycle ESC exhibits a higher radius of oscillation compared to the traditional ESC method. This behavior is due to the control law itself requiring a higher amplitude for the perturbation signal input as observed in \cite{grushkovskaya2025extremum,pokhrel2023higher}. Techniques to attenuate oscillation such as those in \cite{pokhrel2023control, grushkovskaya2024step, yilmaz2025unbiased} can possibly resolve this issue and can be applied to the exponentially convergent unicycle ESC. We also note that attenuation of oscillations in the unicycle-based ESC has been successfully tried in robotic experiments as in \cite{elgohary2025model}.

These results validate the theoretical predictions derived in Section~III and highlight the capability of the proposed approach to handle higher-order objective functions with markedly improved transient response and convergence speed.

\subsection{Experimental Results}
We now extend the validation of the proposed exponentially convergent unicycle ESC design to real-world experiments. The tests are conducted using the same fourth-order objective function employed in the simulation study, with the addition of the high-pass filter (HPF) to improve transient performance. The complete set of parameters for both the traditional and exponentially convergent unicycle ESC designs is listed in Table~\ref{tab:tabexperim}. 


\begin{table}[t]
\centering
\caption{Experimental Parameters for Both Traditional and Exponentially Convergent Unicycle ESC Designs with a Fourth-Order Objective Function}\label{tab:tabexperim}
\begin{tabular}{@{}ll@{}}
\toprule
\textbf{Parameter} & \textbf{Value} \\ 
\midrule
$C_1, \ C_2$ & $1,\ 1$ \\
$a$ & $0.01121$ \\
$c$ & $10$ \\
$\varepsilon$ & $0.2992$ \\
$\Omega$ & $1.4\ \text{rad/s}$ \\
$e$ & $1$ \\
$x(0), \ y(0)$ & $1.6,\ -1.4\ \text{m}$ \\
$x_d, \ y_d$ & $1,\ -2\ \text{m}$ \\
\bottomrule
\end{tabular}
\end{table}

\begin{figure*}[ht]
    \centering
    \begin{subfigure}[b]{0.33\linewidth}
    \centering
    \includegraphics[width=\linewidth]{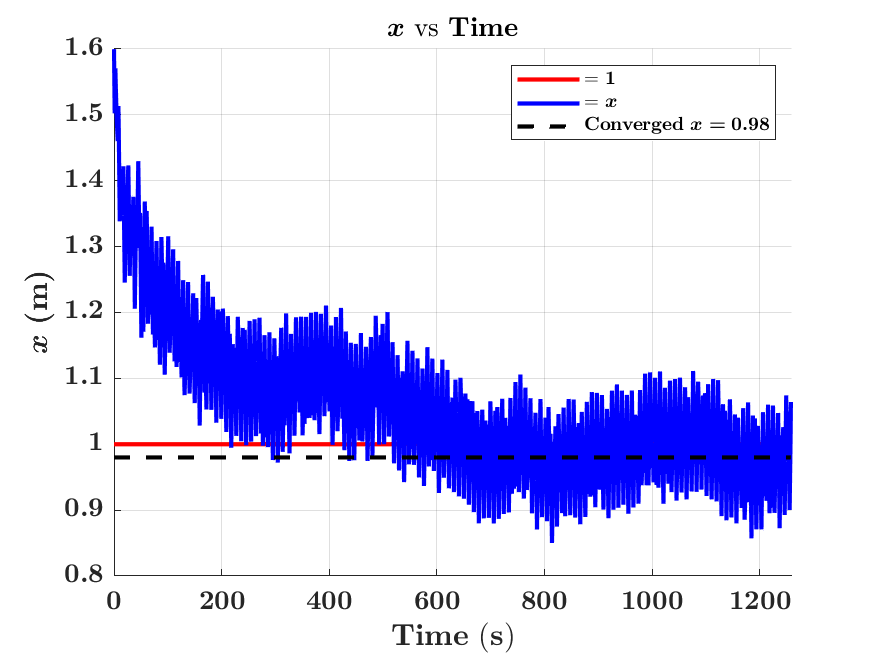}
    \caption{}
    \end{subfigure}\hfill
    \begin{subfigure}[b]{0.33\linewidth}
    \centering
    \includegraphics[width=\linewidth]{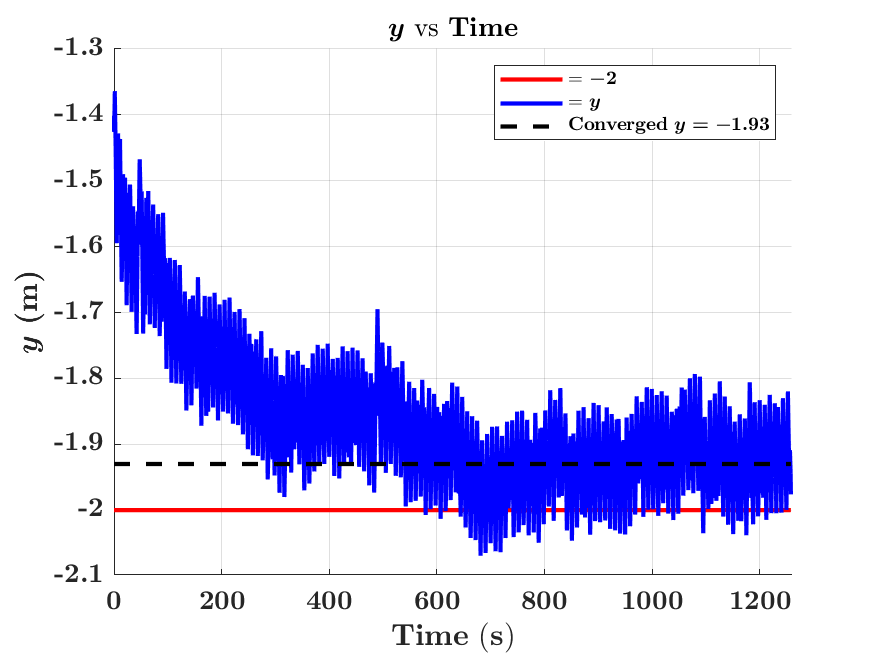}
    \caption{}
    \end{subfigure}\hfill
    \begin{subfigure}[b]{0.33\linewidth}
    \centering
    \includegraphics[width=\linewidth]{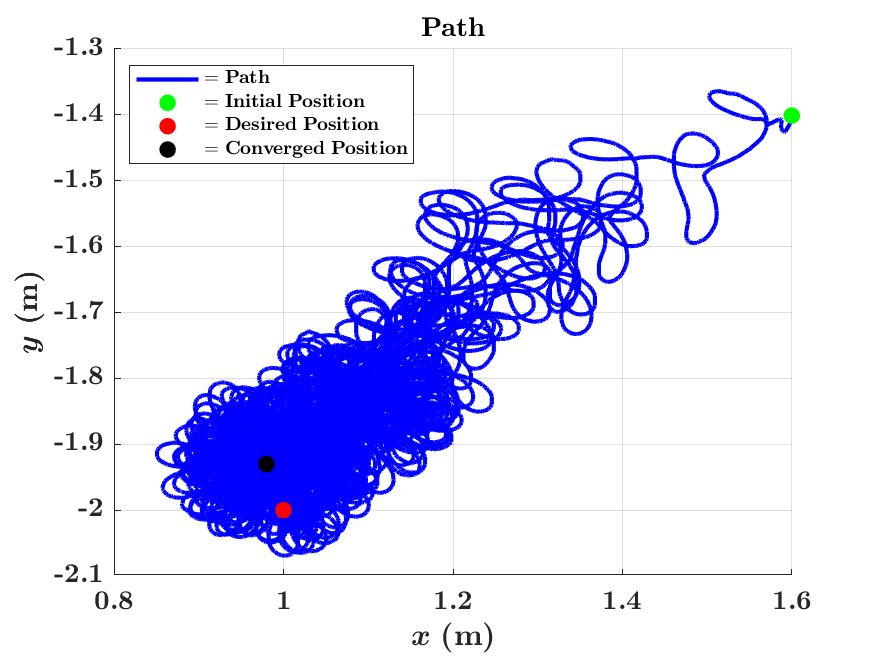}
    \caption{}
    \end{subfigure}
    \caption{Experimental results for the traditional ESC design based on first-order Lie bracket from \cite{durr2013lie} with a fourth-order objective function. (a) $x$ position (b) $y$ position (c) planar trajectory of the robot. Note that the x-position converges to $0.98$ meters (black-dashed), resulting in a 0.02 meter error, and the y-position converges to $-1.93$ (black-dashed), which results in a 0.07 meter error. The true minimum is in red.}
\label{fig:experimental_results_first_ord}
\end{figure*}

\begin{figure*}[ht]
    \centering
    \begin{subfigure}[b]{0.33\linewidth}
    \centering
    \includegraphics[width=\linewidth]{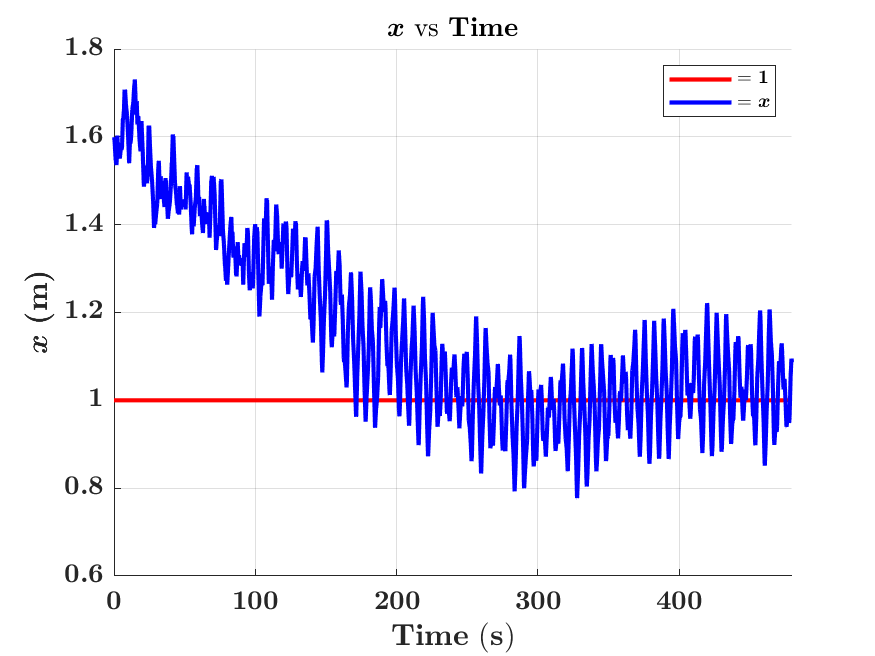}
    \caption{}
    \end{subfigure}\hfill
    \begin{subfigure}[b]{0.33\linewidth}
    \centering
    \includegraphics[width=\linewidth]{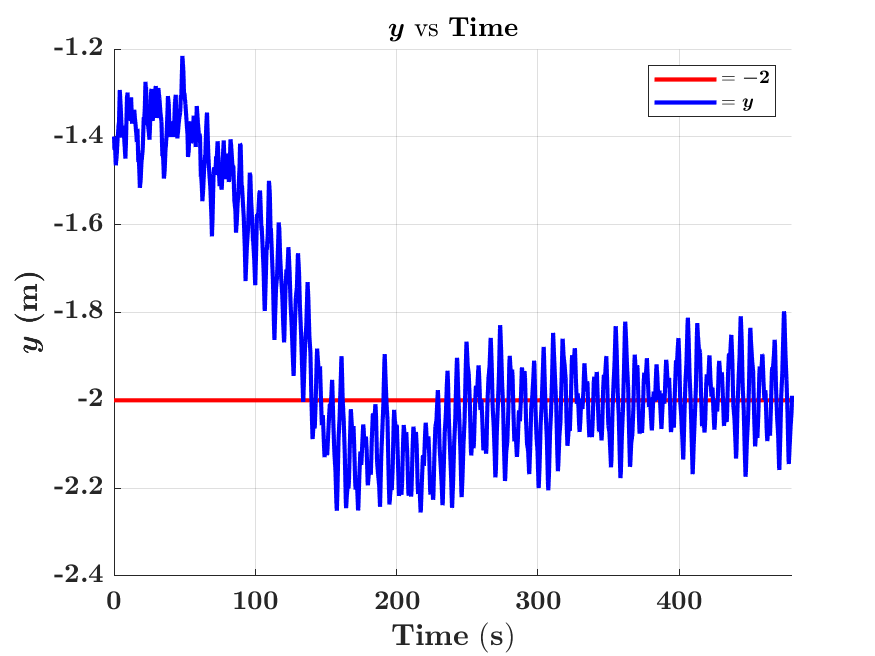}
    \caption{}
    \end{subfigure}\hfill
    \begin{subfigure}[b]{0.33\linewidth}
    \centering
    \includegraphics[width=\linewidth]{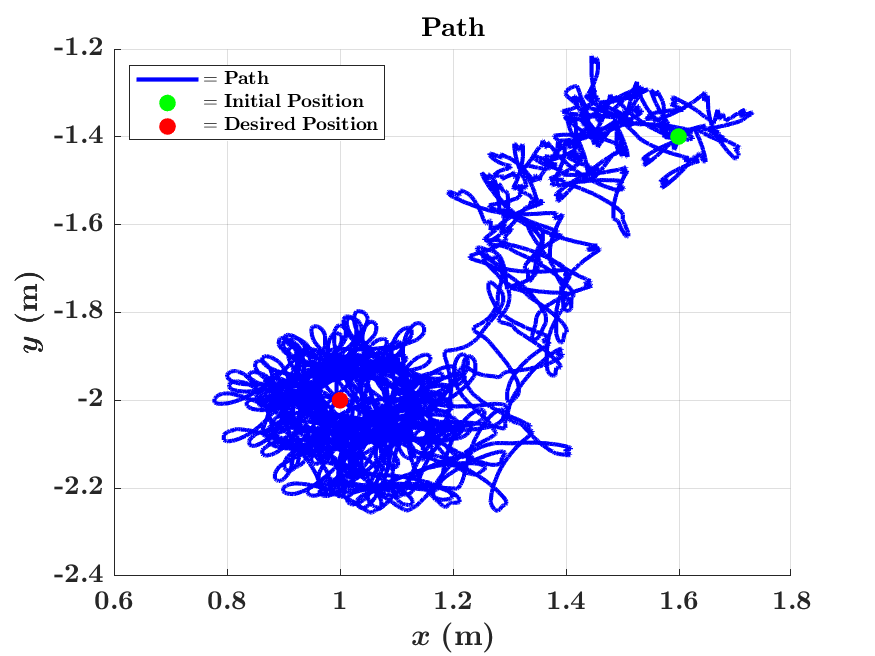}
    \caption{}
    \end{subfigure}\hfill
    \begin{subfigure}[b]{0.33\linewidth}
    \centering
    \includegraphics[width=\linewidth]{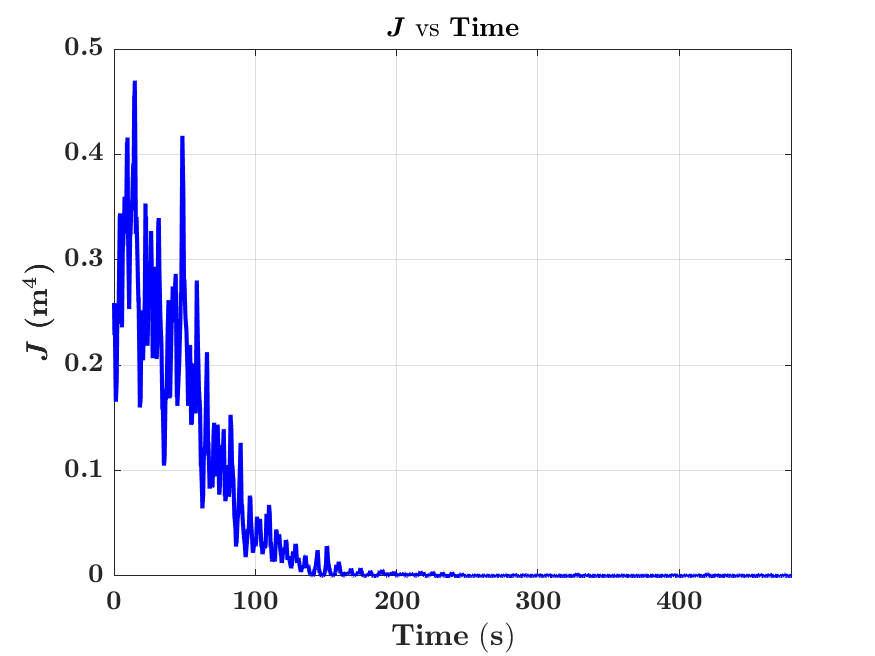}
    \caption{}
    \end{subfigure}\hspace{0.02\linewidth}
    \begin{subfigure}[b]{0.33\linewidth}
    \centering
    \includegraphics[width=\linewidth]{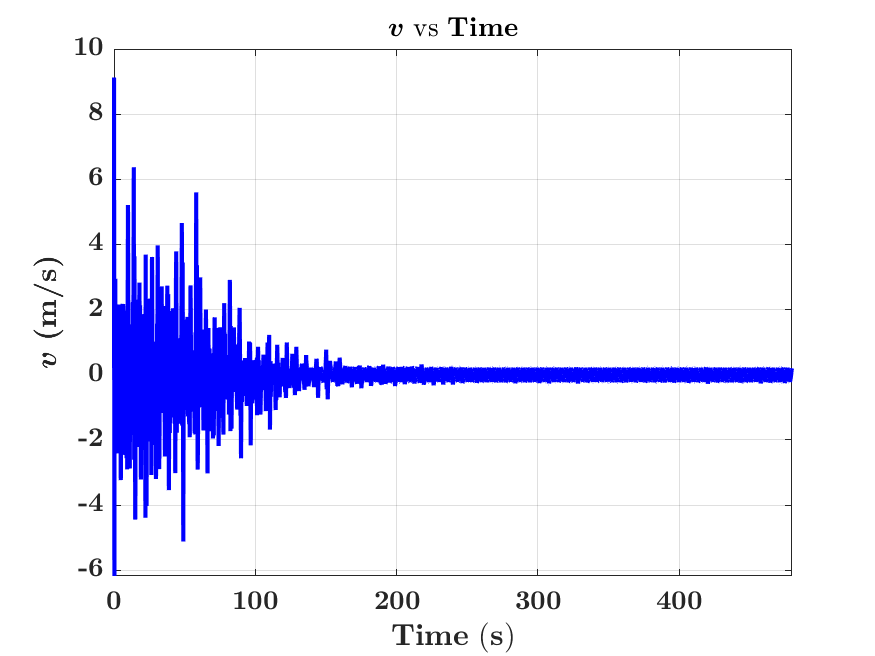}
    \caption{}
    \end{subfigure}

    \caption{Experimental results for the exponentially convergent unicycle ESC design with a fourth-order objective function. (a) $x$ position, (b) $y$ position, (c) planar trajectory of robot, (d) objective function, $J$, (e) linear velocity, $\mathrm{v}$. The reader can watch the experiment at our YouTube channel (\cite{video_known_J}).}
    \label{fig:experimental_results_high_ord}
\end{figure*}

We note here that our design parameters meet the condition of Theorem 1 (see Remark 1). The experimental results are presented in Figures~\ref{fig:experimental_results_first_ord} and~\ref{fig:experimental_results_high_ord}. As shown in Figure~\ref{fig:experimental_results_first_ord}, the traditional ESC method based on first-order Lie bracket from \cite{durr2013lie} fails to converge to the true minimum of the objective function even after 1200~s, instead settling near $(x, y) = (0.98, -1.93)$. Due to the fourth-order objective function, the traditional ESC method will not be able to excite high-order Lie brackets and hence it is not possible to guarantee exponential convergence for this approach. Instead, the traditional ESC method will follow a polynomial decay rate, which will be much slower than an exponential decay rate as discussed and clarified in \cite{grushkovskaya2025extremum}. If given a sufficiently long period of time, then the traditional ESC method in Figure~\ref{fig:experimental_results_first_ord} will eventually converge, as noted in the introduction and discussed in \cite{grushkovskaya2025extremum}. The comparison in Figure~\ref{fig:experimental_results_first_ord} is therefore intended to highlight differences in convergence rates rather than ultimate convergence, with exponential convergence enabling substantially faster practical stabilization within experimentally relevant time horizons.

\begin{figure*}[ht]
    \centering   \includegraphics[width=0.65\linewidth]{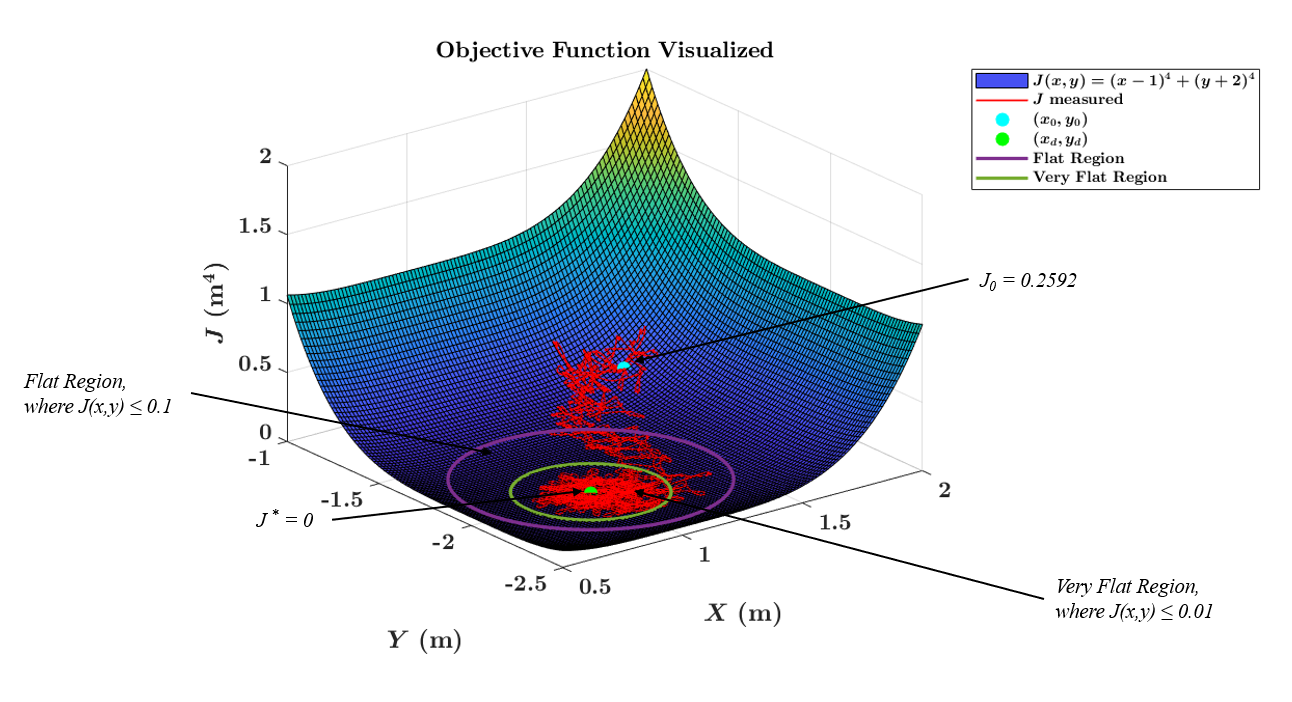}
    \caption{Planar trajectory of the robot with exponentially convergent unicycle ESC design plotted on the objective function distribution. The robot successfully navigates through relatively flat regions where the objective function values are extremely small.}
\label{fig:path_on_objFun_distribution}
\end{figure*}

In contrast, the exponentially convergent unicycle ESC design converges much faster, first reaching the minimum at approximately 350~s and remaining about the true extremum for the duration of the experiment. The new ESC design clearly outperforms the classic literature ESC design and is validated by its ability to reach the minimum in significantly less time than the traditional approach. It is also important to note that a fourth-order objective function exhibits relatively flat regions near the minimum compared to a second-order objective function. The exponentially convergent unicycle ESC design can successfully steer the robot through these flat regions, as evident in Figure~\ref{fig:path_on_objFun_distribution}. Moreover, the variations in the objective function values become very small when the robot is close to the extremum point, confirming that steady-state convergence has been achieved. It is also important to highlight that the proposed design performed well even with expected resolution issues from the motion capturing system feedback due to the very small values of the fourth-order objective function, especially near the extremum. This is a positive indication about the robustness of the proposed design. The reader is directed to watch the experiment in our YouTube channel (\cite{video_known_J}).

\begin{figure*}[t]
    \centering
    \begin{subfigure}[b]{0.33\linewidth}
    \centering
    \includegraphics[width=\linewidth]{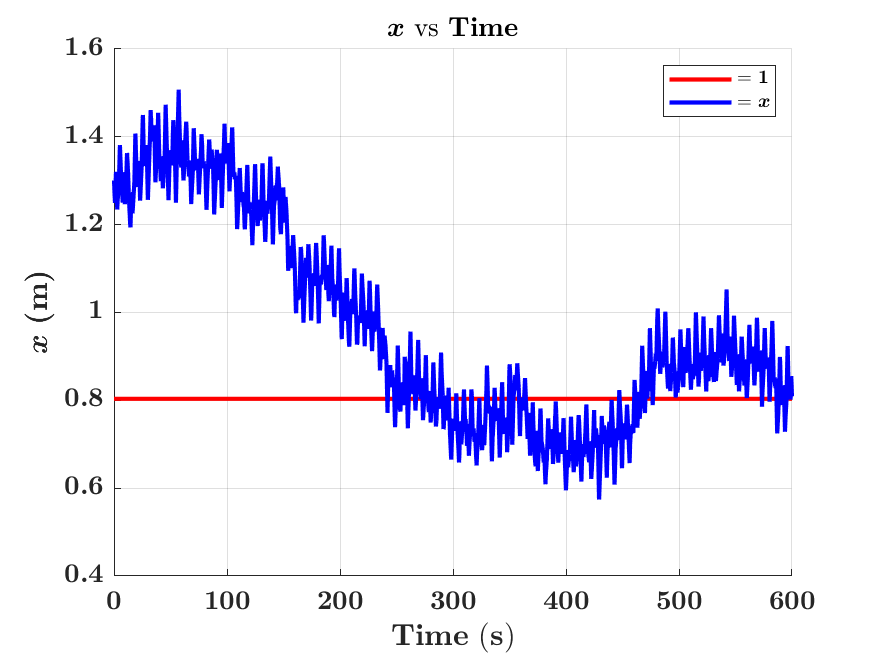}
    \caption{}
    \end{subfigure}\hfill
    \begin{subfigure}[b]{0.33\linewidth}
    \centering
    \includegraphics[width=\linewidth]{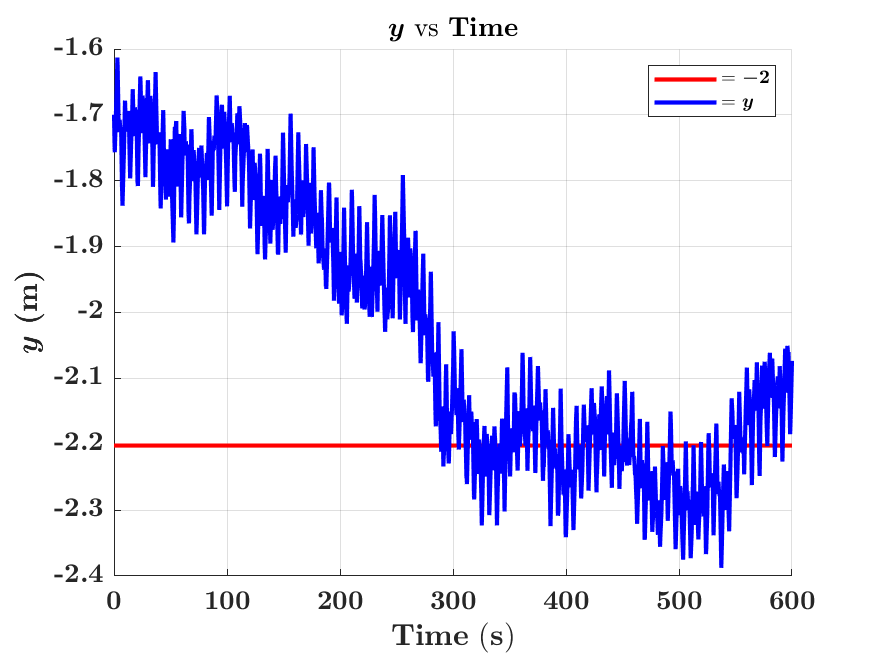}
    \caption{}
    \end{subfigure}\hfill
    \begin{subfigure}[b]{0.33\linewidth}
    \centering
    \includegraphics[width=\linewidth]{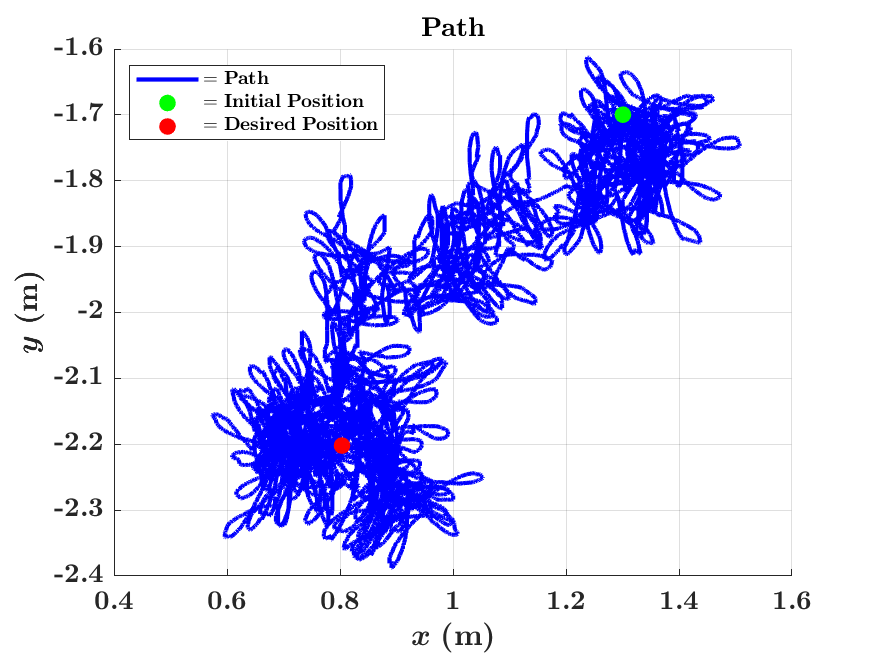}
    \caption{}
    \end{subfigure}\hfill
    \begin{subfigure}[b]{0.33\linewidth}
    \centering
    \includegraphics[width=\linewidth]{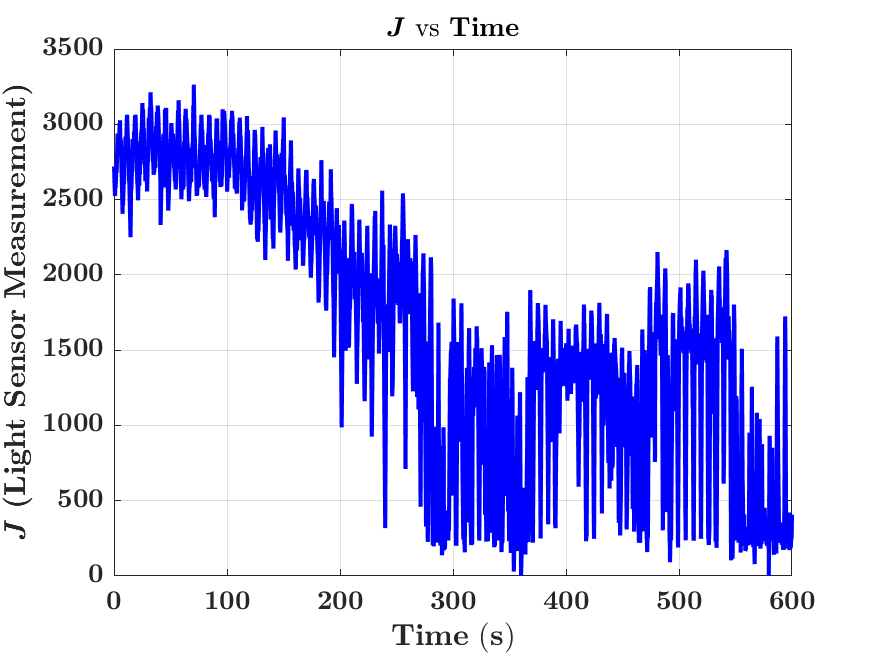}
    \caption{}
    \end{subfigure}\hspace{0.02\linewidth}
    \begin{subfigure}[b]{0.33\linewidth}
    \centering
    \includegraphics[width=\linewidth]{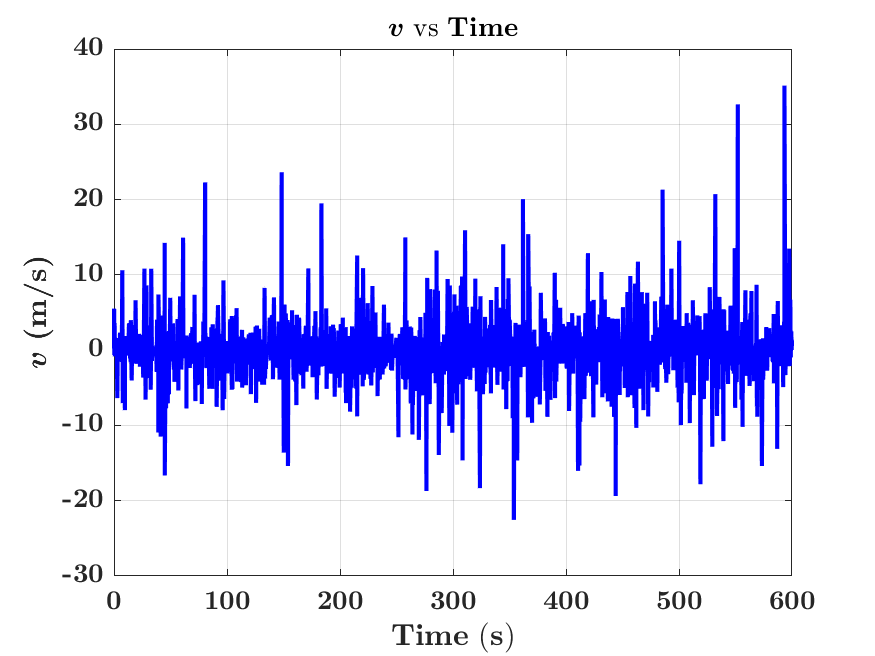}
    \caption{}
    \end{subfigure}

    \caption{Experimental results for light source-seeking with the exponentially convergent ESC model. (a) $x$ position, (b) $y$ position, (c) planar trajectory of robot, (d) objective function, $J$, (e) linear velocity, $\mathrm{v}$. The reader can watch the experiment at our YouTube channel (\cite{video_light}).}
    \label{fig:experimental_results_light}
\end{figure*}

For further verification and validation, we conducted an additional experiment to demonstrate the model-free nature of the proposed exponentially convergent unicycle ESC design. In this test, a light source-seeking task was performed in which the measured light intensity was directly used as feedback to the ESC controller. A high-pass filter (HPF) was again employed to mitigate the influence of natural fluctuations in the light sensor readings. The complete set of parameters used in this experiment is listed in Table~\ref{tab:tablight}. The point of maximum light intensity corresponds to the position where the photoresistive element of the sensor receives the highest illumination, which was determined experimentally to be located at $(x, y) = (0.8035, -2.202)$.
\begin{table}[t]
\centering
\caption{Experimental Parameters for Light Source Seeking with Exponentially Convergent Unicycle ESC Design}\label{tab:tablight}
\begin{tabular}{@{}ll@{}}
\toprule
\textbf{Parameter} & \textbf{Value} \\ 
\midrule
$a$ & $0.006665$ \\
$c$ & $0.001$ \\
$\varepsilon$ & $0.1496$ \\
$\Omega$ & $1.4\ \text{rad/s}$ \\
$e$ & $6$ \\
$x(0), \ y(0)$ & $1.3,\ -1.7\ \text{m}$ \\
\bottomrule
\end{tabular}
\end{table}
The results of this experiment are presented in Figure~\ref{fig:experimental_results_light}. As shown, the robot successfully reaches and oscillates around the location of maximum light intensity. The light sensor measurements decrease as the robot approaches the source, and the mean value of the signal is minimized toward the end of the experiment, despite small oscillations caused by the robot’s continuous motion. These results further confirm that the proposed unicycle ESC design can autonomously steer the robot through an unknown and spatially varying light field, thereby validating its real-time, model-free source-seeking capability. The reader is directed to watch this experiment video in YouTube (\cite{video_light}).
\subsection{Further Discussion on the Selection of $\varepsilon$ in our Experiments}
\label{subsec: discussion_of_eps-selection}
Here, we briefly discuss the selection of $\varepsilon$ in Table~\ref{tab:tabexperim} and Table~\ref{tab:tablight}. In both cases, tuning is performed by selecting a perturbation frequency $\omega$, where $\omega = \frac{2\pi}{\varepsilon}$, and then solving for the corresponding value of $\varepsilon$, as discussed in subsection~\ref{subsec: tuning_control_parms}.

For the case of the known objective function, we selected $\omega = 21$, which yields $\varepsilon = 0.2992$. This choice was guided by experimental observations and prior experience with the robotic platform as in \cite{elgohary2025model,ECC_2024}. While larger values of $\omega$ are preferable in theory, the selected value avoids violating hardware limitations discussed in subsections~\ref{subsec: tuning_control_parms} and~\ref{subsec: summary_of_hardware_limits}. In practice, we verified that higher values such as $\omega = 30$ and $\omega = 42$ were also feasible, but $\omega$ values between 20 and 25 were good enough for the excitation and relatively far from the range of $\omega$ values where hardware limitations are more possible to be encountered.

For the light source-seeking experiment, the objective function is unknown, and therefore the admissible range of $\omega$ (or equivalently $\varepsilon$) may differ from the known objective function case. In particular, the minimum $\omega$ threshold, $\omega^*$, required for effective averaging cannot be guaranteed to be identical for the two systems. That is, the maximum $\varepsilon$ threshold, $\varepsilon^*$, below which $\varepsilon$ must be chosen, cannot be guaranteed to be the same for when the objective function is known compared to when it is unknown since the objective function is the ``output" of the system in ESC designs. Although larger values of $\omega$ remain advantageous in theory (\cite{Maggia2020higherOrderAvg}), we selected $\omega = 42$ to remain within the hardware limitations of the robot. This corresponds to $\varepsilon = 0.1496$.

\subsection{Response to Noise}
In this section, we investigate the robustness of the exponentially convergent unicycle ESC design to noise. To ensure that the numerical simulations capture, as best as possible, the real-world disturbances, we add noise to the measurement of the objective function $J$.

\begin{figure*}[htb]
    \centering
    \begin{subfigure}[b]{0.33\linewidth}
    \centering
    \includegraphics[width=\linewidth]{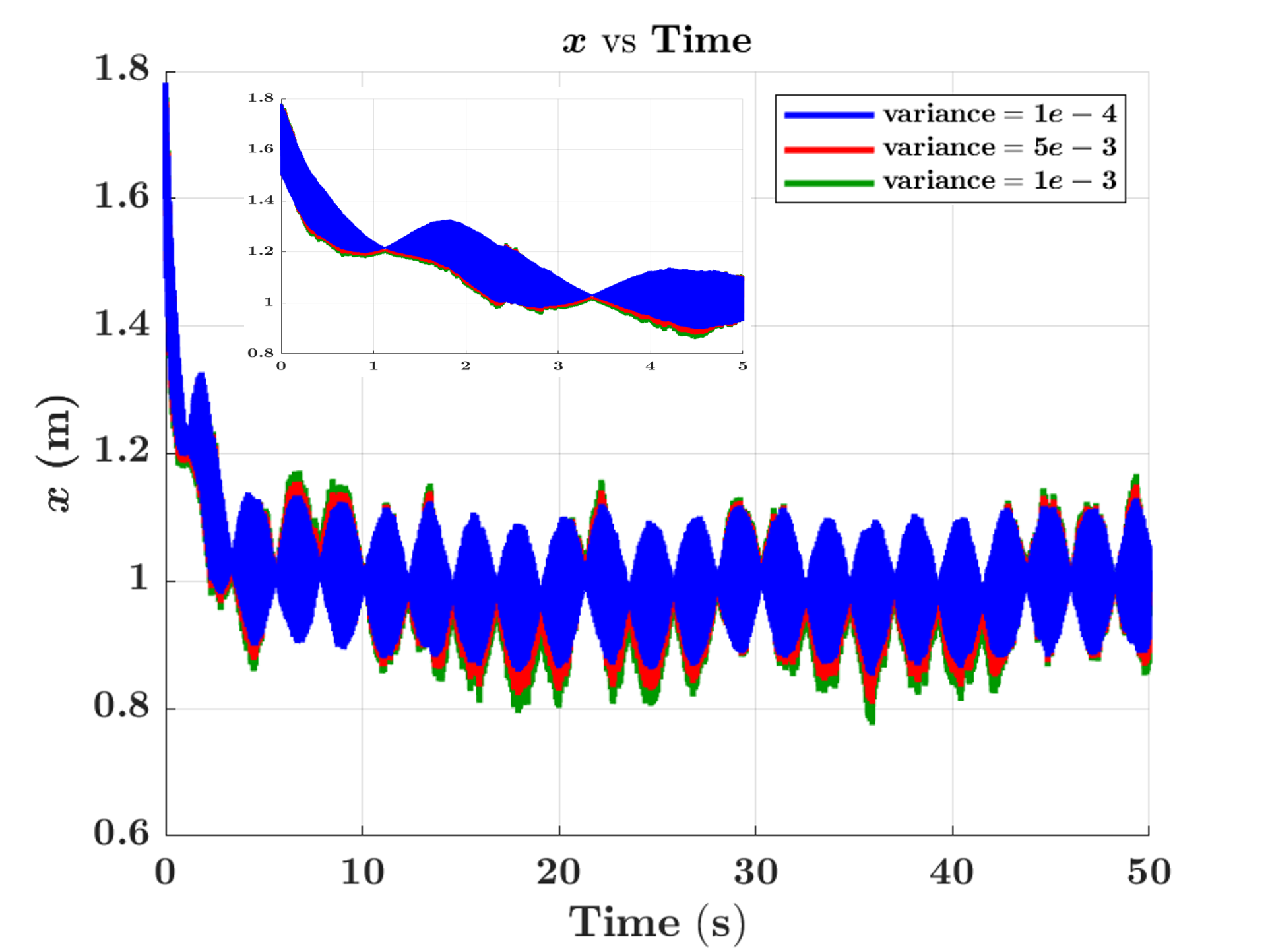}
    \caption{}
    \end{subfigure}\hfill
    \begin{subfigure}[b]{0.33\linewidth}
    \centering
    \includegraphics[width=\linewidth]{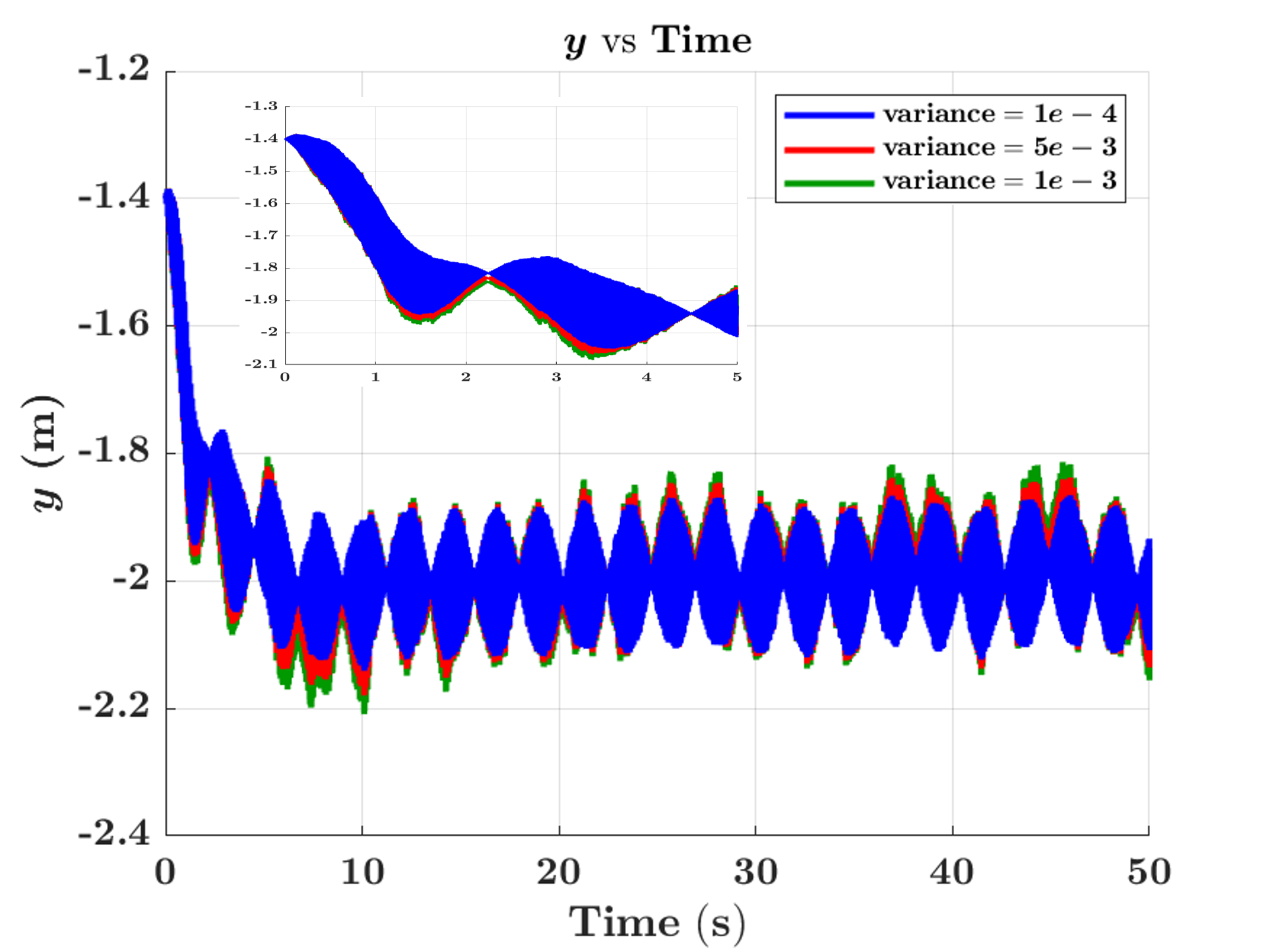}
    \caption{}
    \end{subfigure}\hfill
    \begin{subfigure}[b]{0.33\linewidth}
    \centering
    \includegraphics[width=\linewidth]{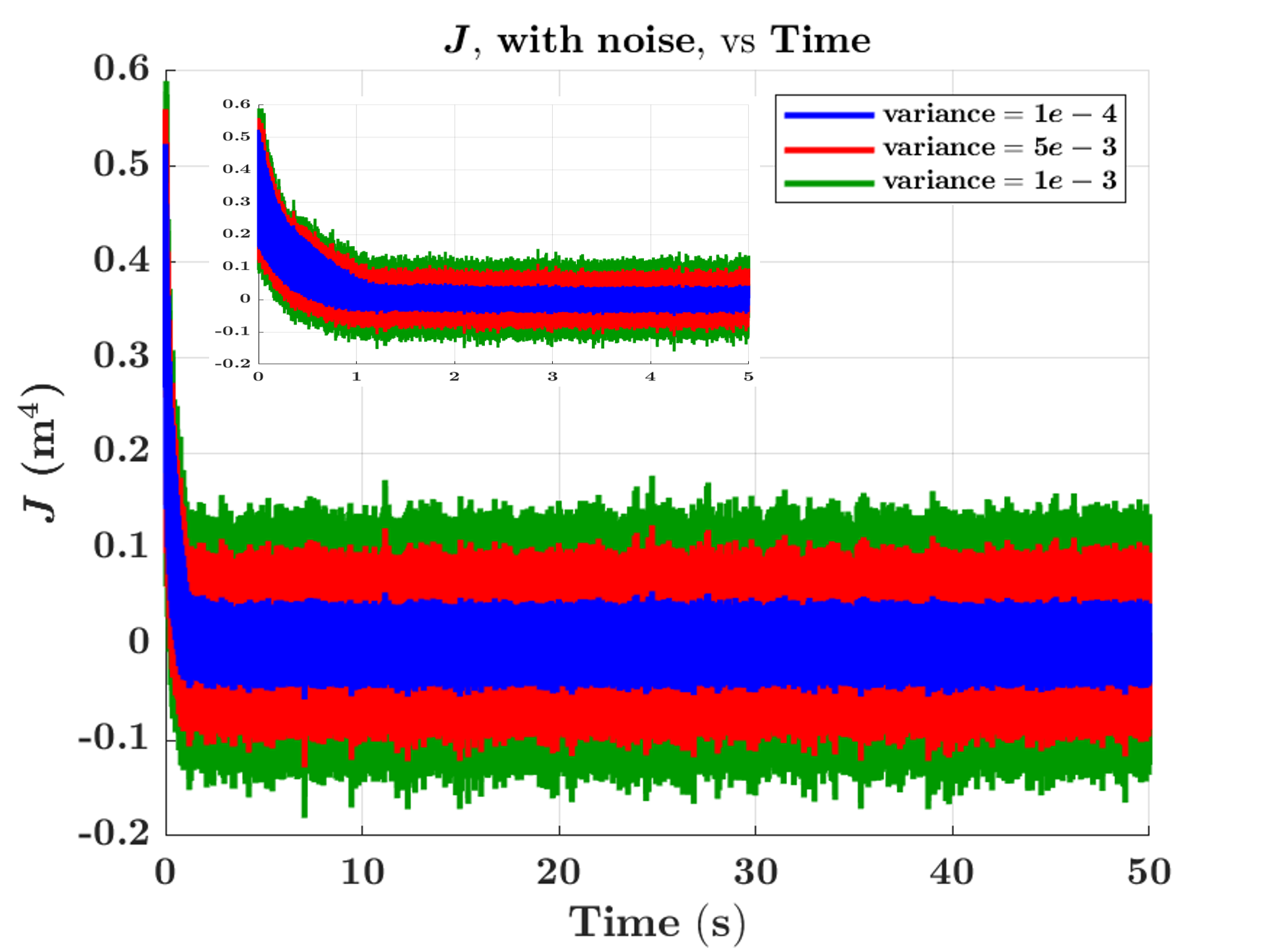}
    \caption{}
    \end{subfigure}
    \caption{Simulation results for the exponentially convergent unicycle ESC with three added measurement noise (low, medium and high). (a) $x$ position, (b) $y$ position, (c) objective function, $J$.}
    \label{fig:simulation_noise_results}
\end{figure*}

We begin by considering numerical simulations for the exponentially convergent unicycle and include three different levels of noise. Noise is added to the objective function using the Simulink \emph{Random Number} block with zero mean and varying variance values to achieve different levels of noise. A fixed seed is used to ensure repeatability when using the \emph{Random Number} block. We conduct numerical simulations using the parameters from Table~\ref{tab:tabsimul} and provide the results in Figure~\ref{fig:simulation_noise_results}. From these results, the unicycle still converges to the desired position even with noise between $10^{-2}$ and $2\times10^{-2}$, which is roughly $38\%$ to $77\%$ of the initial value of the objective function. Moreover, we also report robust response to noise; that is, the lesser the intensity of the noise is, the lesser the radius of oscillation in $x$ and $y$.

\begin{figure*}[htb]
    \centering
    \begin{subfigure}[b]{0.33\linewidth}
    \centering
    \includegraphics[width=\linewidth]{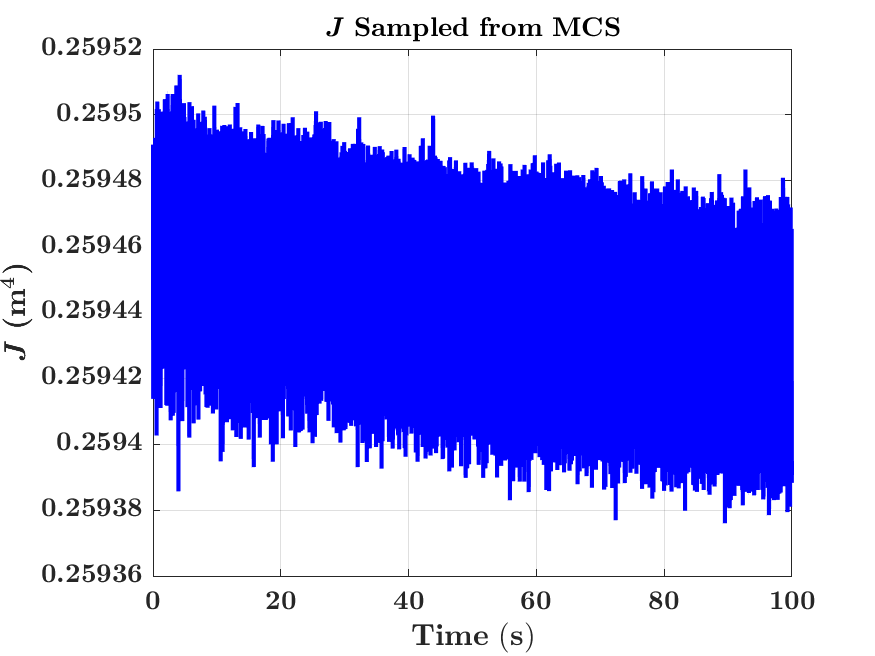}
    \caption{}
    \end{subfigure}\hspace{0.02\linewidth}
    \begin{subfigure}[b]{0.33\linewidth}
    \centering
    \includegraphics[width=\linewidth]{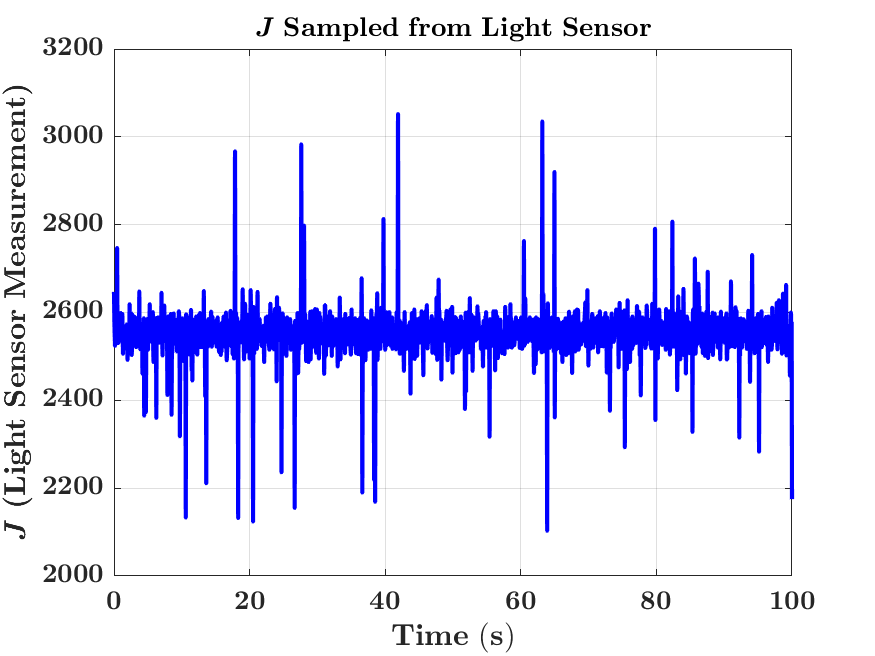}
    \caption{}
    \end{subfigure}\hfill
    \caption{Objective function measurement for the experiment for 100~s with stationary robot at respective initial position in Table~\ref{tab:tabexperim} and Table~\ref{tab:tablight}. (a) Measured objective function from MCS. (b) Measured objective function from the analog light sensor.}
    \label{fig:experiment_noise_results}
\end{figure*}

For the experiment for the known objective function and the light source-seeking, noise was already present in the system from the MCS and analog light sensor. We provide 100~s of objective function measurements while the robot is stationary at the initial position (i.e., without motion) from the MCS and analog light sensor in Figure~\ref{fig:experiment_noise_results}. Yet, as reported in the previous subsection, all experiments were successful in convergence despite of these inherent noises.

\subsection{Response to Delay}
We investigate the robustness of the exponentially convergent unicycle ESC design to delay. To best relate the numerical simulations to real-world applications, we include a delay in the objective function, $J$.

\begin{figure*}[htb]
    \centering
    \begin{subfigure}[b]{0.33\linewidth}
    \centering
    \includegraphics[width=\linewidth]{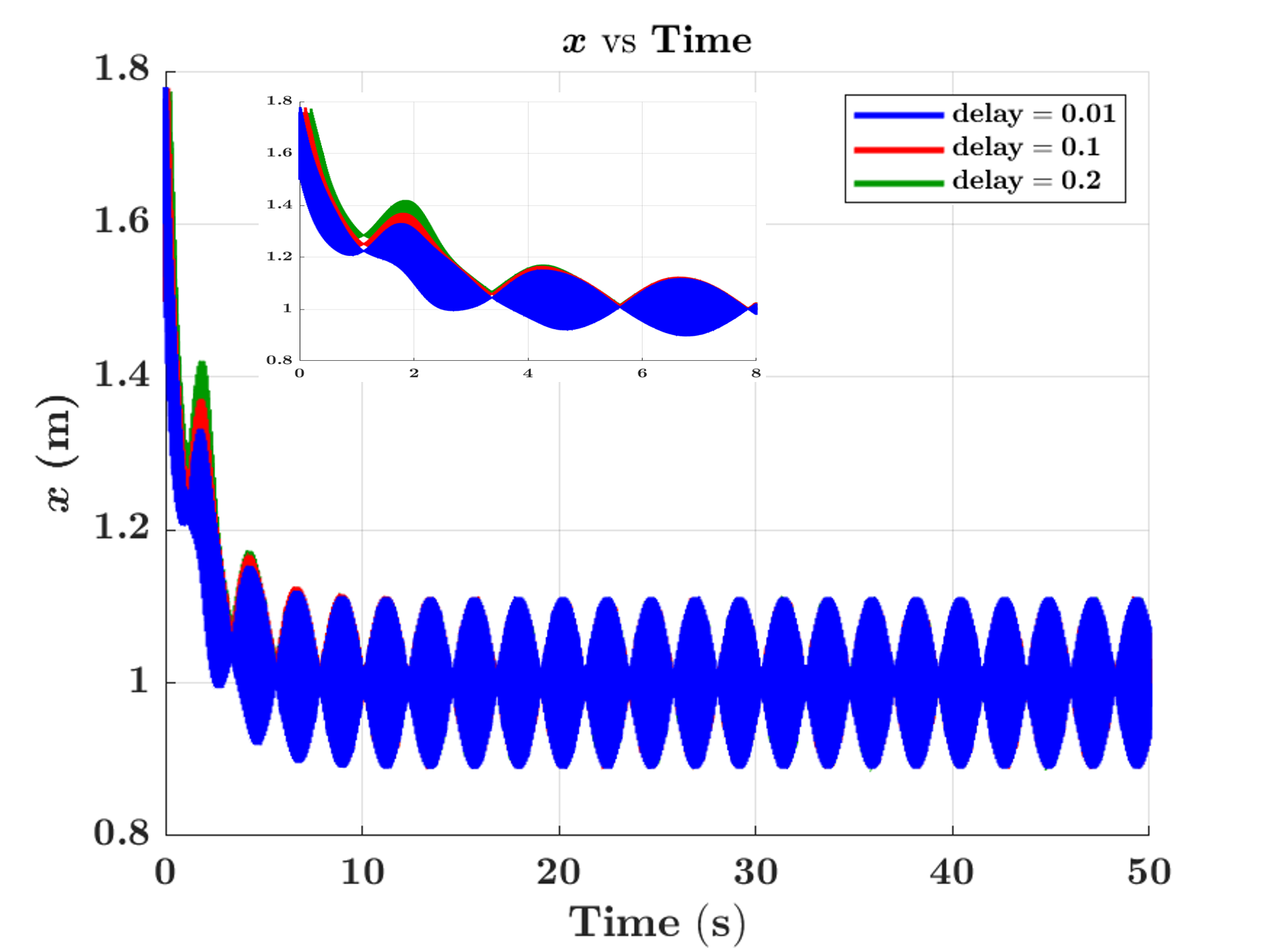}
    \caption{}
    \end{subfigure}\hfill
    \begin{subfigure}[b]{0.33\linewidth}
    \centering
    \includegraphics[width=\linewidth]{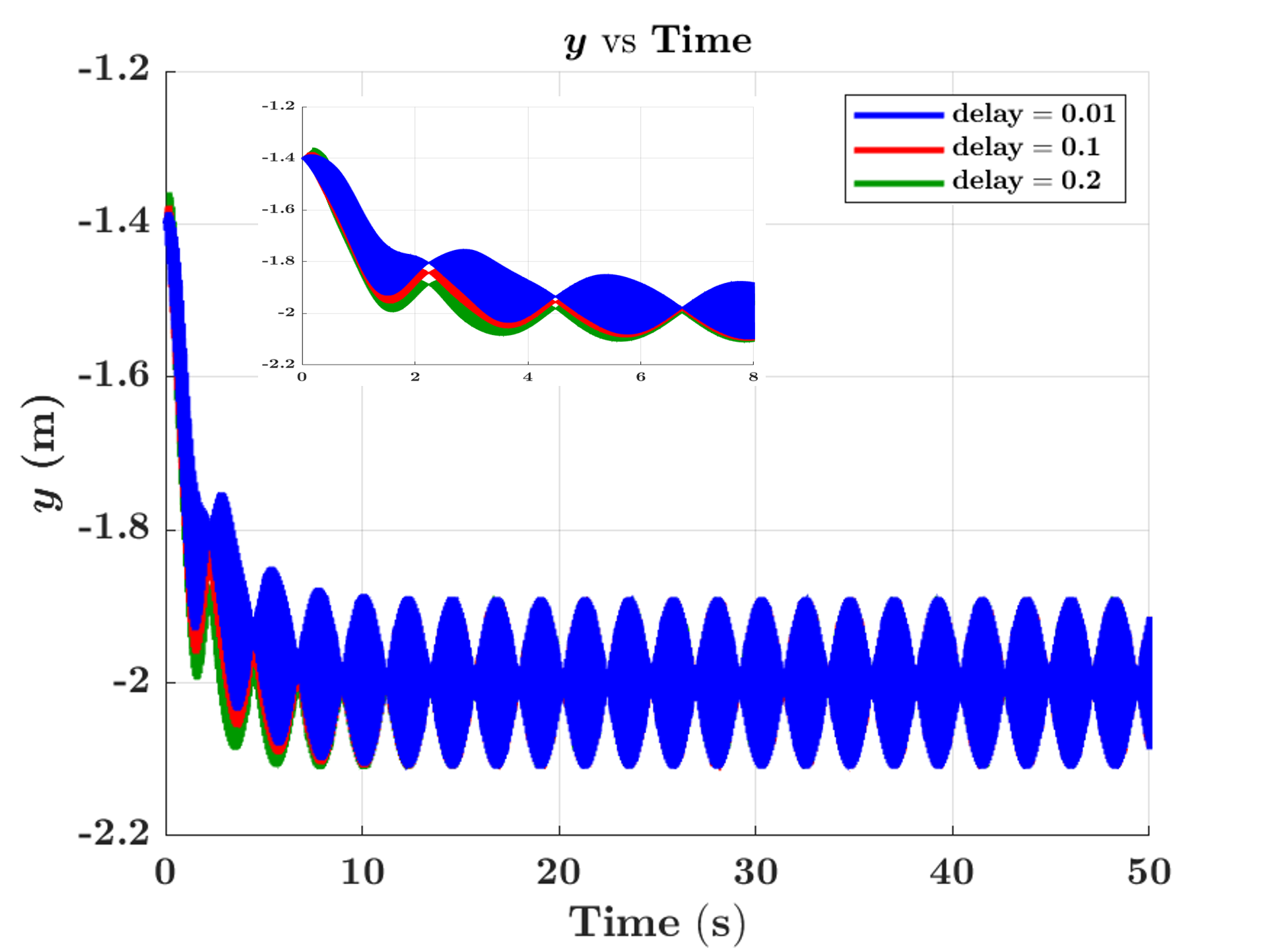}
    \caption{}
    \end{subfigure}\hfill
    \begin{subfigure}[b]{0.33\linewidth}
    \centering
    \includegraphics[width=\linewidth]{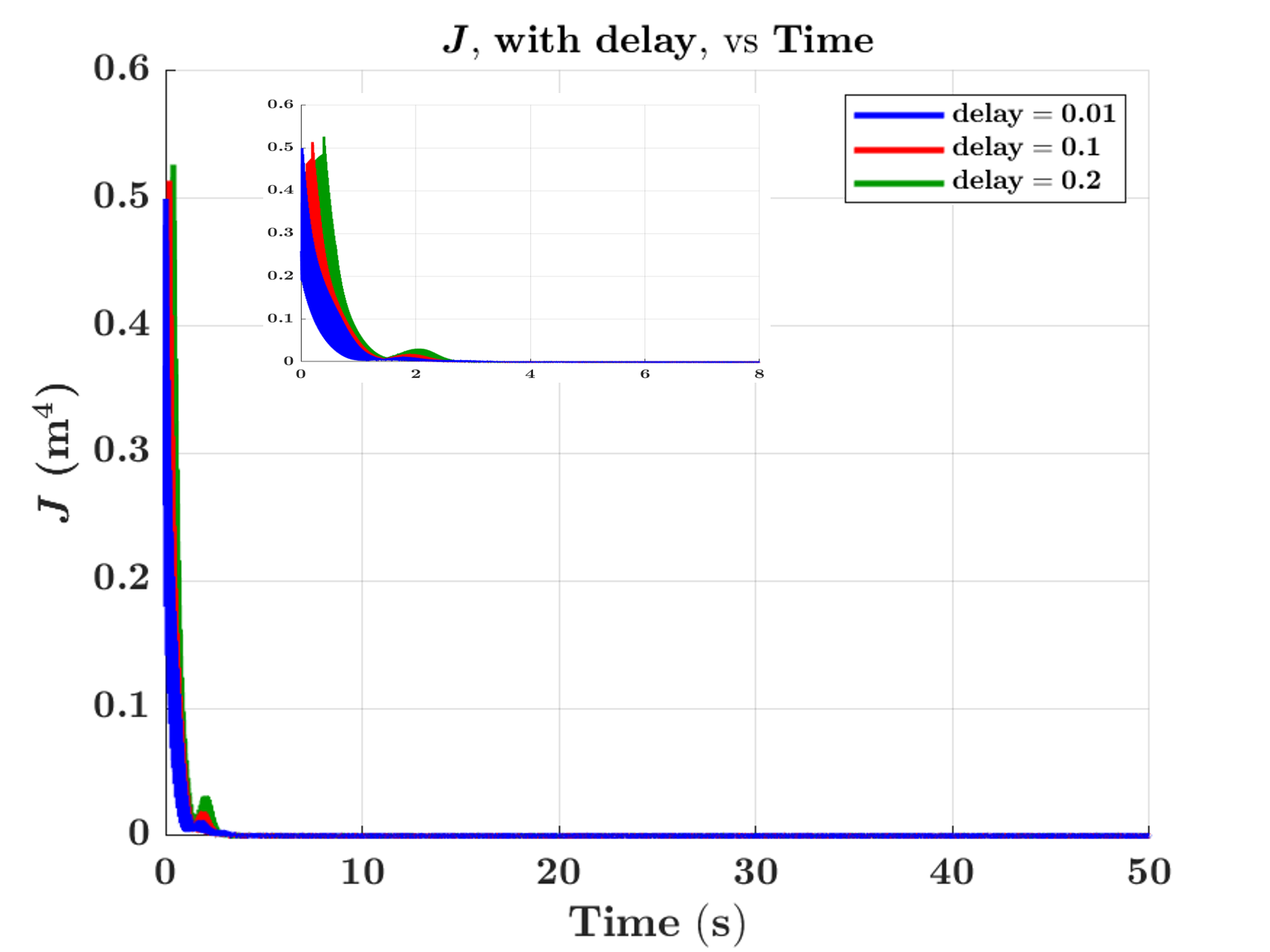}
    \caption{}
    \end{subfigure}
    \caption{Simulation results for exponentially convergent unicycle ESC with delay. (a) $x$ position, (b) $y$ position, (c) objective function, $J$.}
    \label{fig:simulation_delay_results}
\end{figure*}

We consider three different levels of delay for the numerical simulations. The Simulink \emph{Transport Delay} block, with an initial value equivalent to the initial value of the objective function, is used. We use the parameters from Table~\ref{tab:tabsimul} and present the results in Figure~\ref{fig:simulation_delay_results}. We observe that even with a delay of up to $0.2$~s, the exponentially convergent unicycle still converges to the desired position quickly.

Furthermore, we note that we did not compensate for delay in the experiment that is inherently present in the system due to the latency between the obtained measurements and the latency in setting the commanded velocity by the robot. As seen in Figure~\ref{fig:experimental_results_high_ord} and Figure~\ref{fig:experimental_results_light}, the robot was able to converge to the desired position regardless of the delays present in the system.

\subsection{Response to Different Initial Values in Simulations and Experiment}
We investigate the response of the exponentially convergent unicycle ESC with different initial conditions to prove the robustness of the system, and draw some observations about potential gap between simulations and experiments.

\begin{figure*}[htb]
    \centering
    \begin{subfigure}[b]{0.33\linewidth}
    \centering
    \includegraphics[width=\linewidth]{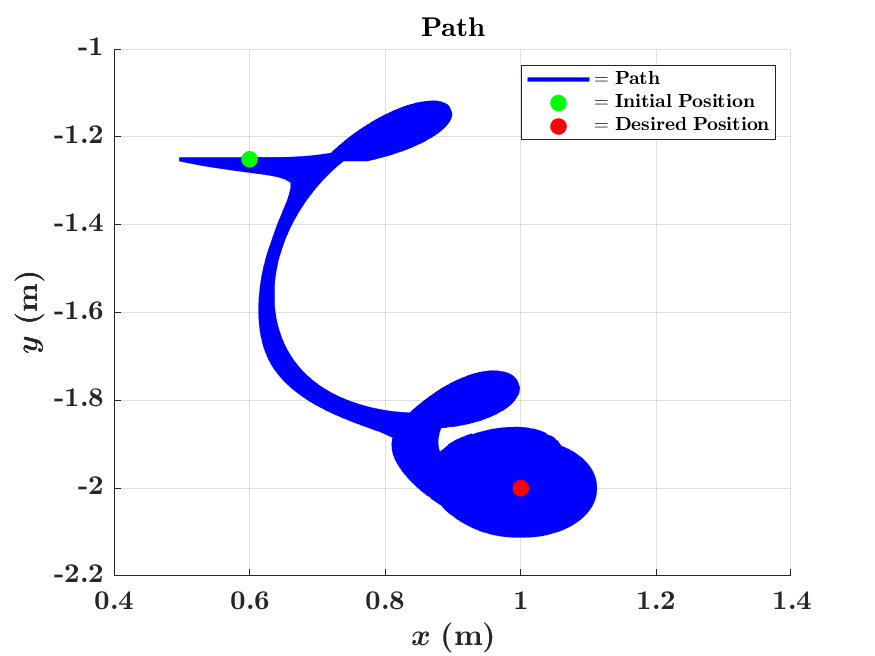}
    \caption{}
    \end{subfigure}\hfill
    \begin{subfigure}[b]{0.33\linewidth}
    \centering
    \includegraphics[width=\linewidth]{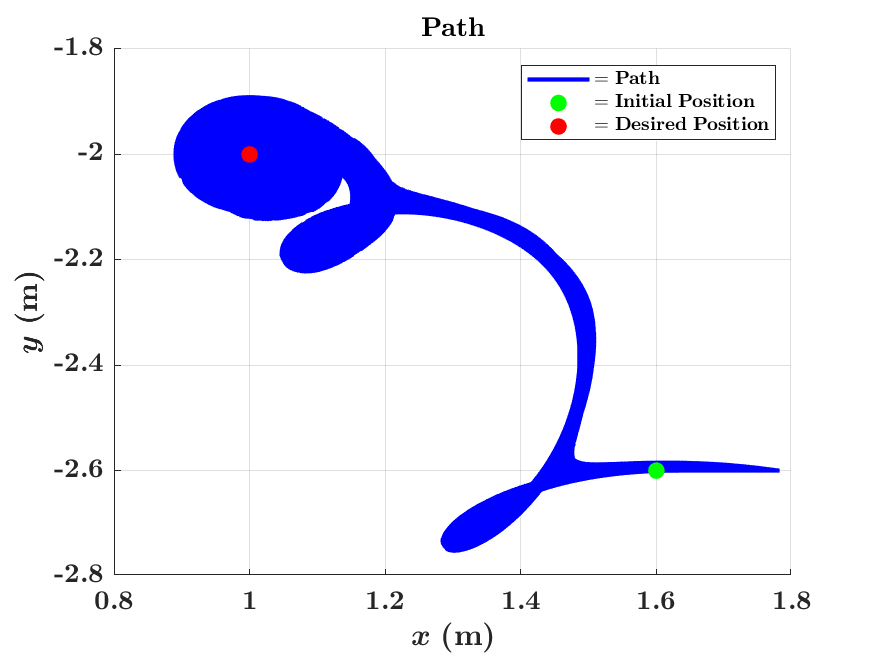}
    \caption{}
    \end{subfigure}\hfill
    \begin{subfigure}[b]{0.33\linewidth}
    \centering
    \includegraphics[width=\linewidth]{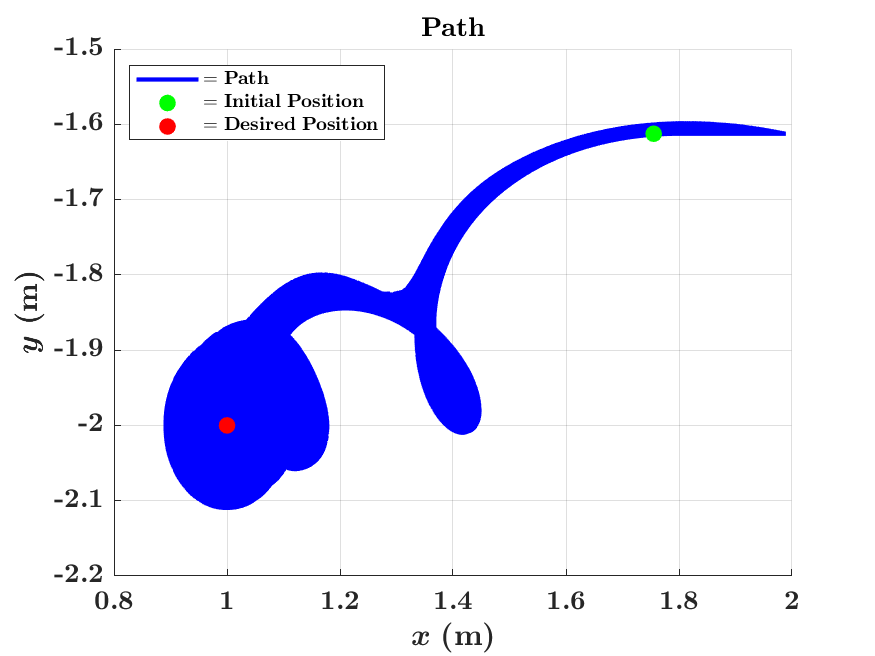}
    \caption{}
    \end{subfigure}
    \caption{Simulation results for different initial conditions representing the same radius from the extremum to test consistency in convergence of the exponentially convergent unicycle ESC with respect to direction. (a) $x(0),\ y(0) = 0.6, \ -1.251$, (b) $x(0),\ y(0) = 1.6, \ -2.6$, (c) $x(0),\ y(0) = 1.7552, \ -1.6126$.}
    \label{fig:simulation_initcond_results}
\end{figure*}

\begin{figure*}[htb]
    \centering
    \begin{subfigure}[b]{0.33\linewidth}
    \centering
    \includegraphics[width=\linewidth]{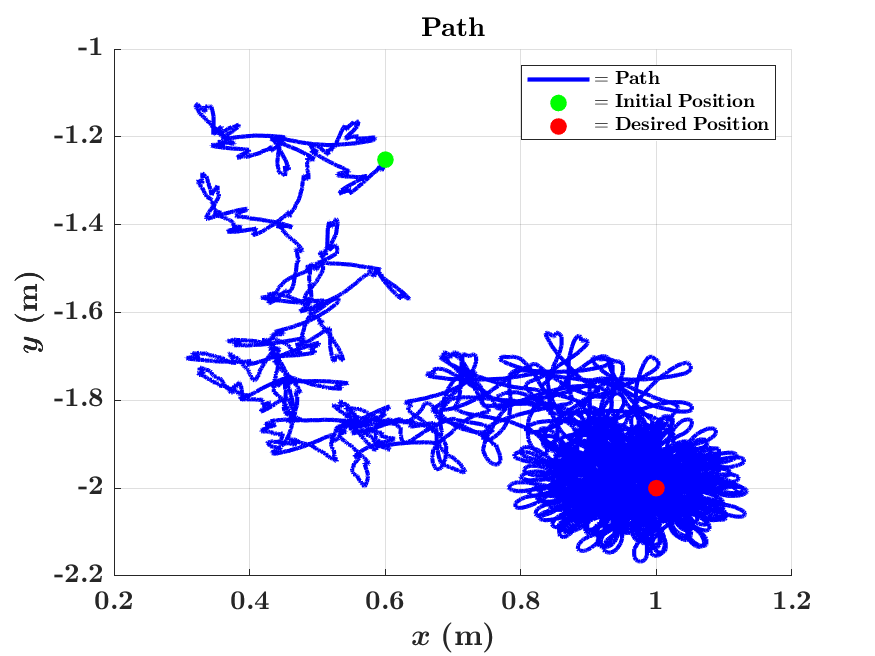}
    \caption{}
    \end{subfigure}\hfill
    \begin{subfigure}[b]{0.33\linewidth}
    \centering
    \includegraphics[width=\linewidth]{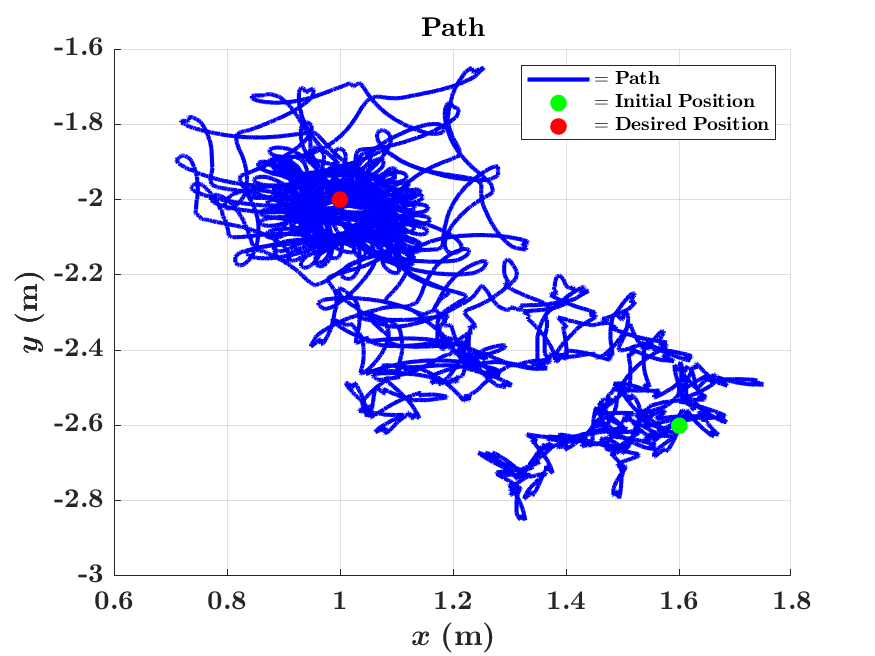}
    \caption{}
    \end{subfigure}\hfill
    \begin{subfigure}[b]{0.33\linewidth}
    \centering
    \includegraphics[width=\linewidth]{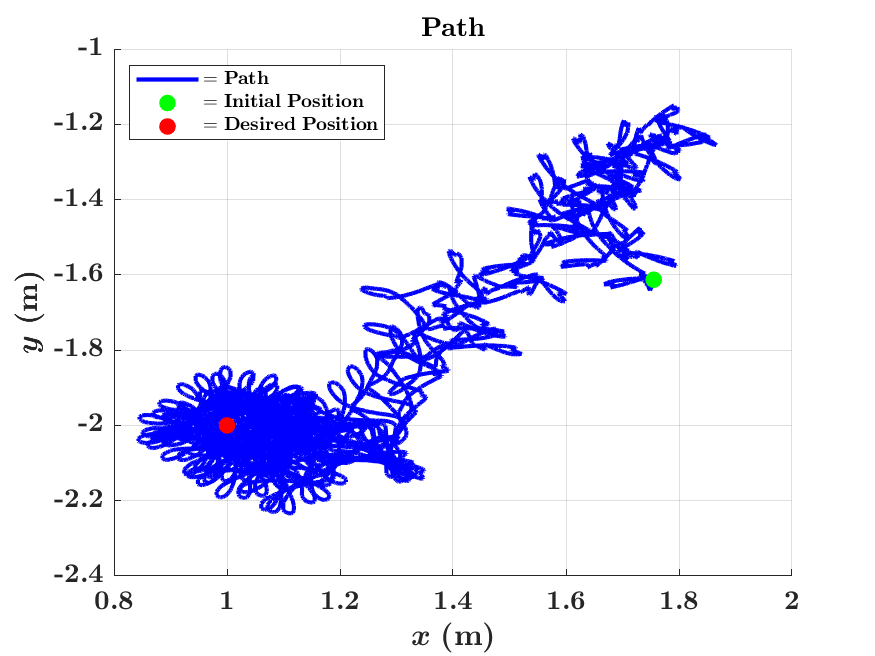}
    \caption{}
    \end{subfigure}
    \caption{Experimental results for different initial conditions representing the same radius from the extremum to test consistency in convergence of the exponentially convergent unicycle ESC with respect to direction. (a) $x(0),\ y(0) = 0.6, \ -1.251$, (b) $x(0),\ y(0) = 1.6, \ -2.6$, (c) $x(0),\ y(0) = 1.7552, \ -1.6126$.}
    \label{fig:experiment_initcond_results}
\end{figure*}

We begin by testing three different starting positions in both simulation and experiment. First, we consider the following three sets of initial conditions: (a) $x(0),\ y(0) = 0.6, \ -1.251$, (b) $x(0),\ y(0) = 1.6, \ -2.6$, (c) $x(0),\ y(0) = 1.7552, \ -1.6126$. These initial conditions were selected to reflect different positions that are about the same radius from the extremum point to test directional consistency in convergence. The simulations are conducted using the parameters from Table~\ref{tab:tabsimul}, and the experiments are conducted using the parameters from Table~\ref{tab:tabexperim}. The simulation results are provided in Figure~\ref{fig:simulation_initcond_results}, and the experimental results are provided in Figure~\ref{fig:experiment_initcond_results}. In both simulation and experiments, we show that regardless of the initial position direction, the unicycle or robot can successfully converge to the desired position, demonstrating consistency and robustness of convergence from different directions.

\begin{figure*}[ht]
    \centering   \includegraphics[width=0.5\linewidth]{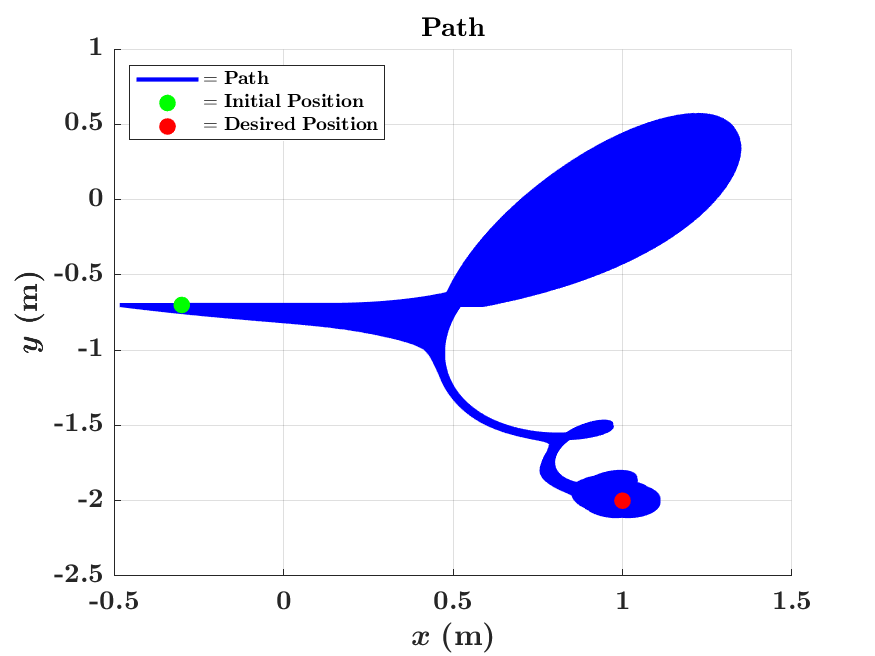}
    \caption{Simulation result for an initial position that is farther away from the desired position with the exponentially convergent unicycle ESC. $x(0),\ y(0) = -0.3, \ -0.7$.}
\label{fig:simulation_initcond_faraway}
\end{figure*}

Furthermore, we extend the distance between the initial and desired positions and show, in simulation, the ability of the exponentially convergent unicycle ESC to reach the desired position from a farther distance. The initial condition for this simulation is taken as, $x(0),\ y(0) = -0.3, \ -0.7$ and the result is shown in Figure~\ref{fig:simulation_initcond_faraway}. The unicycle is still able to converge to the desired position even when starting farther away from the desired position.

In the experiment, due to the hardware limitations of the TurtleBot3 robot, we were not able to test initial conditions that were farther away from the desired position. Specifically, the maximum linear velocity of the TurtleBot3 limits the admissible amplitude of the velocity perturbation signal. When the objective function value is large (i.e., far from the extremum), the commanded velocity saturates, which suppresses the effective excitation needed to realize the corresponding Lie bracket direction. As a result, initial conditions that are too far from the desired position cannot be reliably explored in experiment.

Overall, the exponentially convergent unicycle ESC does not show dependence on the initial position directions with respect to the extremum in both simulations and experiments. However, a gap between theory/simulations on the one hand and experiments on the other, was observed as the exponentially convergent unicycle has hardware limitations preventing it from converging if the initial condition is set very far from the extremum; a solution for that could be a future work that develops a bounded update rate design upgrade for the exponentially convergent ESC, analogous to the bounded update design upgrade in \cite{scheinker2017bounded} with respect to the traditional ESC design in \cite{durr2013lie}.

\subsection{Comparison between Simulation and Experiment Results}
We further investigate a direct comparison between the simulation results and the experimental results to evaluate if the performance of the exponentially convergent unicycle ESC matches the theoretically predicted exponential convergence rate provided in Section 3 that is illustrated by simulations.

\begin{figure*}[ht!]
    \centering
    \begin{subfigure}[b]{0.33\linewidth}
    \centering
    \includegraphics[width=\linewidth]{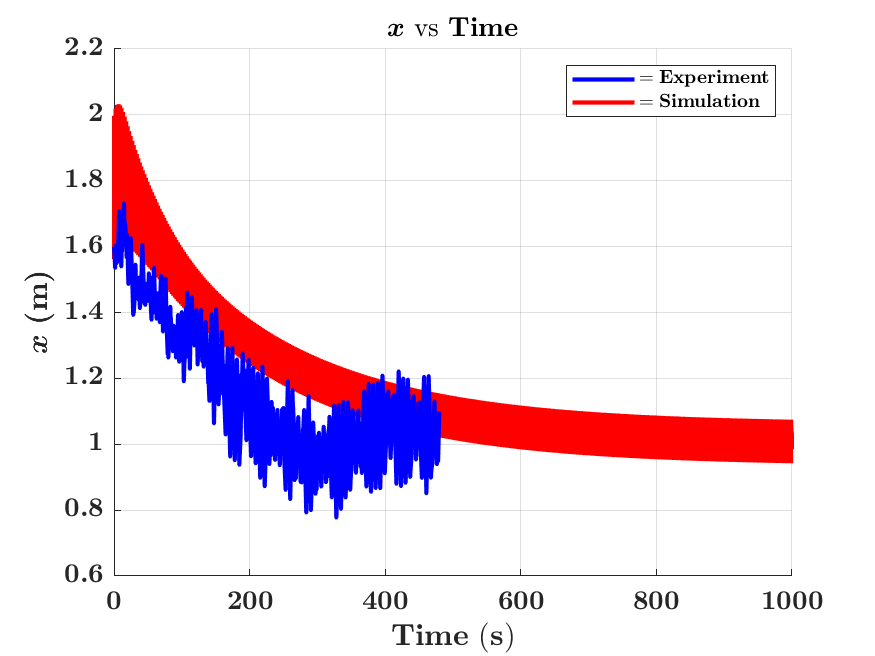}
    \caption{}
    \end{subfigure}\hfill
    \begin{subfigure}[b]{0.33\linewidth}
    \centering
    \includegraphics[width=\linewidth]{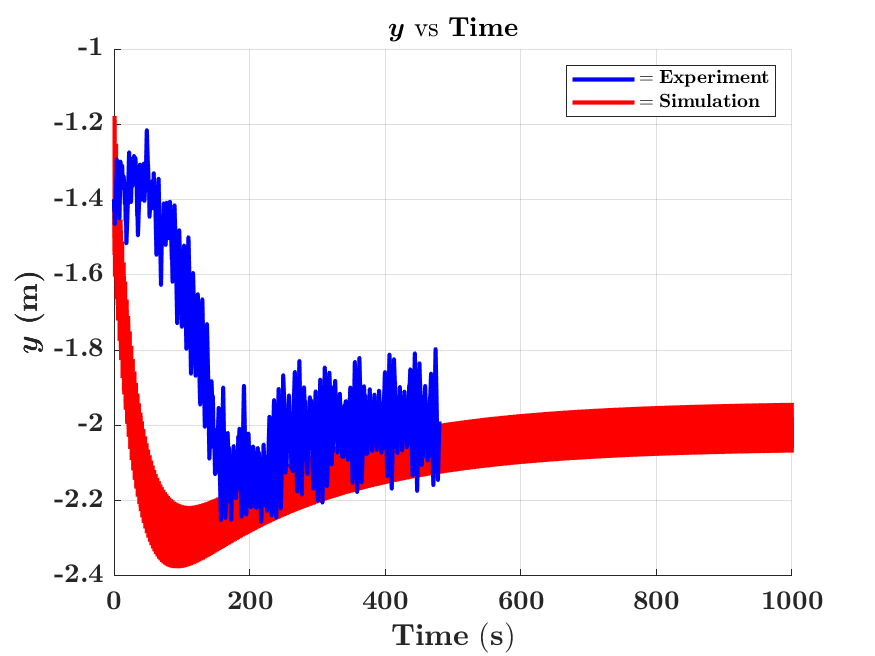}
    \caption{}
    \end{subfigure}\hfill
    \begin{subfigure}[b]{0.33\linewidth}
    \centering
    \includegraphics[width=\linewidth]{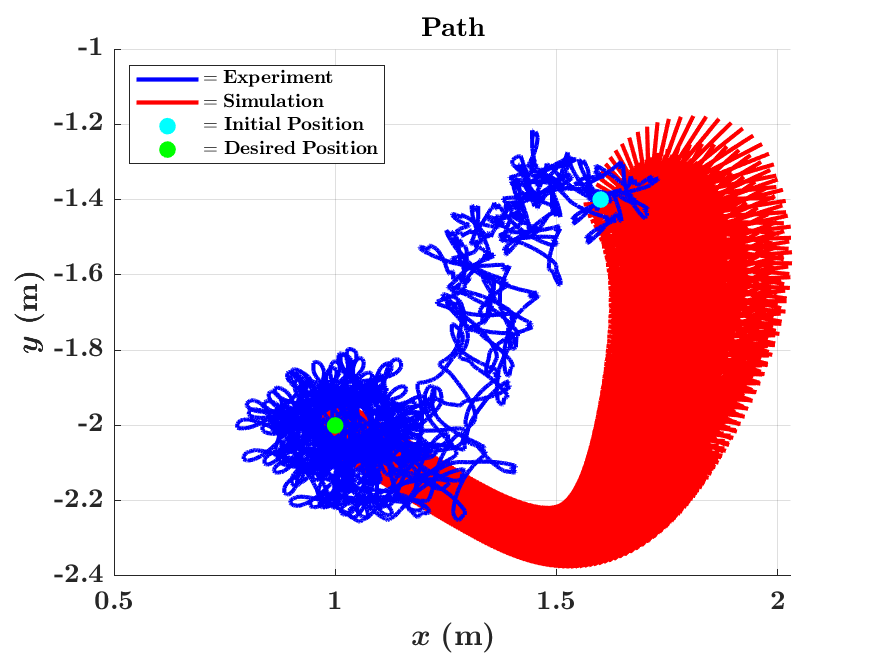}
    \caption{}
    \end{subfigure}
    \caption{Comparison between experiment-based system (blue) and simulation-based system (red) with close, but still higher initial convergence rate demonstrates that the proposed exponentially convergent unicycle design follows the predicted theoretical behavior.}
    \label{fig:sim_exp_comp_results}
\end{figure*}

To begin, we note that the unicycle dynamics model does not perfectly capture the true dynamics of the TurtleBot3 robot. Inertial terms are not captured in the unicycle dynamics as the unicycle is treated as a point mass. Also, dissipation terms, such as friction, are not captured within the unicycle dynamics. For these reasons, it is not possible to directly equate the simulation and experiment parameters. That is, Table~\ref{tab:tabsimul} and Table~\ref{tab:tabexperim} are not aligned in terms of the parameters, partly due to the differences in the unicycle model and the real-world system of TurtleBot3. To compare the simulation results with its theoretically predicted behavior vs. the experimental result, we aimed to match the initial convergence rate for the states between the two systems (simulation-based and experiment-based) in a similar fashion to what could be found in ESC literature when one compares the convergence rate of two different ESC systems, see for example the comparisons in \cite{elgohary2025extremum}. Analyzing Table~\ref{tab:tabsimul} and Table~\ref{tab:tabexperim}, we see that at $t = 0$~s, $\dot{y}(0)$ is the same for both systems and the only difference between the two systems is in $\dot{x}(0)$ with the simulation-based system possessing a much higher initial convergence rate in $\dot{x}$ compared to the experiment-based system. That is, for the experiment-based system, $\dot{x} = 0.2199$, and for the simulation-based system, $\dot{x}(0) = 705.724$. We now introduce a new set of parameters, in Table~\ref{tab:tabsimul_lowconv}, to reduce the initial convergence rate of $\dot{x}$ in the simulation-based system as much as possible to allow for a relatively fairer comparison with the experiment-based system. The experimental parameters are the same as in Table~\ref{tab:tabexperim}. With the new parameters in Table~\ref{tab:tabsimul_lowconv} for the simulation-based system, the initial convergence rate for $\dot{x}$ becomes $\dot{x}(0) = 5.5382$. Now, the results that compare the experiment-based system and the simulation-based system, which, theoretically, is guaranteed to have an exponential convergence rate, are provided in Figure~\ref{fig:sim_exp_comp_results}. We observe that the simulation-based system, even with relatively higher initial convergence rate in $\dot{x}$ (i.e., more advantageous initial speed), is slower than the experimental-based system. This shows that the proposed unicycle ESC follows the theoretical predictions in Section 3 and their simulations; however, one need to keep in mind the reported limitation in Subsection 4.6 where the real-world system may not converge for initial conditions that are very far from the extremum.

\begin{table}[t]
\centering
\caption{Simulation Parameters for Lower $\dot{x}$ Convergence Rate}\label{tab:tabsimul_lowconv}
\begin{tabular}{@{}ll@{}}
\toprule
\textbf{Parameter} & \textbf{Value} \\ 
\midrule
$C_1, \ C_2$ & $1,\ 1$ \\
$a$ & $0.1$ \\
$c$ & $0.5$ \\
$\varepsilon$ & $0.075$ \\
$\Omega$ & $1.4\ \text{rad/s}$ \\
$x(0), \ y(0)$ & $1.6,\ -1.4\ \text{m}$ \\
$x_d, \ y_d$ & $1,\ -2\ \text{m}$ \\
\bottomrule
\end{tabular}
\end{table}

We also note that for the model-free light source-seeking experiment, we use different parameters compared to the known objective function experiment and simulations. This is in part because the light distribution is not known, and the scale of the measurements received from the light sensor differs from the cases when the objective function is known as fourth-order polynomial, hence differences in tuning are expected.

\subsection{Summary of Observed Hardware Limitations}
\label{subsec: summary_of_hardware_limits}
We note that two primary hardware limitations arise from our experiments.
\begin{itemize}
    \item \textbf{Minimum value of $\varepsilon$.} Theoretically, smaller values of $\varepsilon$ (equivalently, larger perturbation frequencies $\omega$) are preferable. In practice, however, there is a lower bound on $\varepsilon$ imposed by the physical limitations of the robot. Large values of $\omega$ challenge the ability of the motors to track the commanded velocity input. Despite this limitation, our experiments demonstrate successful convergence for values of $\varepsilon$ that respect the physical frequency constraints of the platform, as shown in Figures~\ref{fig:experimental_results_high_ord}, \ref{fig:experimental_results_light}, and \ref{fig:experiment_initcond_results}. Additional details on the selection of $\varepsilon$ are provided in subsection~\ref{subsec: discussion_of_eps-selection}. Moreover, since the control loop is executed with a fixed sampling time, the Nyquist frequency also imposes a lower bound on admissible values of $\varepsilon$.

    \item \textbf{Maximum linear velocity of the robot.} The maximum linear velocity of the TurtleBot3 imposes a physical saturation limit on the velocity perturbation signal. As a result, the experimental basin of attraction differs from the theoretical one. This effect becomes particularly evident when the initial condition is far from the extremum. While the unicycle model in simulation can converge from distant initial positions (Fig.~\ref{fig:simulation_initcond_faraway}), the physical robot cannot replicate this behavior due to velocity saturation. For fourth-order objective functions, the objective value grows rapidly with distance from the minimum, leading to large commanded velocities that saturate at the robot’s physical limit. In this regime, the effective excitation required to realize the Lie bracket direction is suppressed, and the robot may oscillate or drift. In contrast, for the experimental cases shown in Figures~\ref{fig:experimental_results_high_ord} and~\ref{fig:experiment_initcond_results}, the initial objective values remain sufficiently small to allow the perturbation signal to excite the correct descent direction. Consequently, the initial conditions in Fig.~\ref{fig:simulation_initcond_faraway} cannot be replicated experimentally using the current platform.
\end{itemize}

\section{Conclusion and future work}
This paper provided a novel, first-of-its-kind exponentially convergent unicycle design inspired by recent results on higher-order Lie bracket approximations (\cite{grushkovskaya2025extremum,pokhrel2023higher}). We proved exponential stability of the proposed design for objective functions that behave locally as a fourth-degree polynomial. The theoretical results were validated through extensive simulations and, more importantly, through real-world robotic experiments. In simulations, we demonstrated robust performance of the proposed design under different initial conditions as well as in the presence of measurement noise and delays. In experiments, we validated not only the ability of the proposed design to operate with known objective functions and varying initial conditions, but also its capability to operate in a fully model-free setting through a light source-seeking experiment. Furthermore, a case study comparing experimental and simulation results was presented to confirm that the observed experimental behavior follows the theoretically predicted exponential convergence trends. These results demonstrate that exponentially convergent ESC laws based on higher-order Lie bracket approximations can be practically realized on physical robotic platforms, despite hardware limitations and measurement imperfections.

In future work, we aim to generalize the proposed design to objective functions with higher polynomial degrees and to develop modified ESC structures that reduce oscillations and enforce bounded update rates while preserving exponential convergence properties.

\bibliography{ifacconf}

\end{document}